\documentclass[12pt]{article}

\usepackage{graphicx}
\usepackage{amssymb}
\usepackage{icomma} 
\usepackage{caption}
\usepackage{tikz}
\usetikzlibrary{matrix}
\usepackage{mathrsfs}
\usepackage{empheq}
\usepackage{amsmath}
\usepackage[utf8]{inputenc} 
\usepackage[T1]{fontenc}    
\usepackage{hyperref}       
\usepackage{url}            
\usepackage{booktabs}       
\usepackage{amsfonts}       
\usepackage{nicefrac}       
\usepackage{microtype}      
\usepackage{lipsum}

\newtheorem{theorem}{Theorem}[section]

\newtheorem{definition}[theorem]{Definition}

\newtheorem{proposition}[theorem]{Proposition}

\title{Existence of a Variational Principle for PDEs with Symmetries and Current Conservation}

\author{Markus Dafinger}
\date{\today}

\begin{document}
\maketitle

\begin{abstract}
We prove that under certain assumptions a partial differential equation can be derived from a variational principle. It is well-known from Noether\rq{}s theorem that symmetries of a variational functional lead to conservation laws of the corresponding Euler-Lagrange equation. We reverse this statement and prove that a differential equation which satisfies sufficiently many symmetries and corresponding conservation laws leads to a variational functional whose Euler-Lagrange equation is the given differential equation.
\end{abstract}


\section{Introduction}
\label{intro}
In this paper we prove that a differential equation which satisfies certain kinds of symmetries and conservation laws can be written as an Euler-Lagrange equation, that is, it can be derived from a variational principle. Before we make this statement more precise we provide some motivation and explain possible applications.\\
\hspace*{0.5cm} The variational principle is very important in theoretical physics and it seems that all the differential equations which describe the laws of nature on a fundamental level can be derived from a variational principle, as for example, Einstein\rq s field equations and the equations in the standard model. Feynman\rq s path integral formalism also relies on it which can be used to quantize a classical system. Since Lagrange formulated this principle in 1788 \cite{bib49}, it is hard to imagine doing theoretical physics without it.\\ 
\hspace*{0.5cm} Because the variational principle has been used very successfully in physics, it is desirable to have some explanation why it works so well, and to understand whether we can state it as a kind of superordinate axiom in physics, or whether we can prove it as a consequence of other conditions, like symmetries and conservation laws. Symmetries, like Lorentz-, Poincaré-, or certain kinds of gauge invariance can be possibly assumed for any fundamental differential equation in physics. For the conservation laws, we may assume energy-, mass-, momentum-, angular momentum-, or charge conservation, among others.\\ 
\hspace*{0.5cm} The above question, whether a differential equation which satisfies certain symmetries and conservation laws is necessarily derivable from a variational principle, was first formulated by Floris Takens in 1977. He also proved three different theorems in which he answered this question affirmatively \cite{bib17}. Others have generalized his proofs in a series of papers which was initiated by I. M. Anderson and J. Pohjanpelto in 1994. Scalar PDEs of second order are considered in \cite{bib12}, system of PDEs (SPDEs) of first order in \cite{bib15}, SPDEs of second order in \cite{bib14,bib16}, SPDEs of third order in \cite{bib1}, and SPDEs of arbitrary order, but only polynomial expressions, are investigated in \cite{bib13}. In the above mentioned articles there are usually additional assumptions in their statements, beside the assumptions of order, scalar PDEs, SPDEs, and polynomial structure. In \cite{bib15,bib14,bib16,bib1} very strong restrictions for the symmetries and conservation laws are made. For the symmetries, very often some kind of translation-, and gauge invariance is assumed. Since there is a connection between symmetries and conservation laws (similar as in Noether\rq s theorem), the conservation law assumptions in these theorems are also of a very special type. For example, gauge invariance leads to the condition that the differential equation must be divergence-free. In this article we generalize the case of second order SPDEs in the way that we do not restrict to very special symmetries, especially as it is done for second order SPDEs in \cite{bib14,bib16}. Symmetries are describe by a set of vector fields, and it can be shown that such a set is a Lie algebra \cite[p.177]{bib36}. Our main theorem delivers a proof for any non-abelian, or abelian Lie algebra of symmetries which spans a certain tangent space pointwise. It generalizes and unifies previously known results by Takens, Anderson, Pohjanpelto, Manno, and Vitolo of a similar nature. A more detailed list of our generalizations can be found below after we have introduce some more notation.\\
\hspace*{0.5cm} We also want to emphasize here that our approach to Takens\rq{} problem relies on a kind of Cartan\rq{}s formula which is, in our opinion, essential. This formula, see \eqref{lack3}, can be understood very easily when the differential equation is transformed to a weak formulation for certain kinds of test functions. Surprisingly, this kind of Cartan\rq{}s formula, and especially the weak formulation, is not directly used in many other papers which consider Takens\rq{} problem and it makes it much more challenging to understand the nature of the problem. Although the connection to Cartan\rq{}s formula is mentioned sometimes, the explicit generalized exterior derivatives $D_{n+1}$ and $D_{n+2}$ in the Lie derivative $\mathcal{L}_V=D_{n+1}\lrcorner V+V\lrcorner D_{n+2}$ are not stated, and it turns out that they are not obvious to find in the variational sequence for differential forms, since they do not coincide with the Euler-Lagrange operator $E_{n+1}$ and Helmholtz operators $E_{n+2}$, see \eqref{euler} and \eqref{helm}.\\
\hspace*{0.5cm} A differential equation $f(x,u(x),Du(x),...)=0$ contains independent coordinates $x$ and dependent coordinates $u$. The concept of independent- and dependent coordinates is realized by a fiber bundle $\pi:E\to M$, where $x$ are local coordinates on the base manifold $M$ and $u$ are local coordinates of the fibers of $E$. We let $n=\operatorname{dim}M$ and $n+m=\operatorname{dim}E$. Derivatives $Du$, $D^2u$ and so on are then local coordinates on certain fibers of the jet bundle $J^kE$. Symmetries are $\pi$-projectable vector fields on $E$ which are lifted to vector fields on $J^kE$ and this lift is called jet prolongation.\\ 
\hspace*{0.5cm} As already mentioned, instead of investigating the differential equation $f=0$ directly, we consider a weak formulation of the differential equation, that is, $\int f\varphi dx=0$, where $\varphi$ are test functions with certain properties, that is, they correspond to $\pi$-vertical vector fields on $E$, as we will see below. The rough formulation of our main theorem is the following (the precise formulation can be found in Theorem \ref{t47}):\\
\begin{theorem}\label{th0}
Let $\pi:E\to M$ be a fiber bundle with base dimension $n$ and fiber dimension $m$, where $n,m\in\mathbb{N}$ are arbitrary, and $M$ is oriented. Furthermore, let 
\begin{align*}
K_f(u;\varphi)=\int_M f_\alpha(x,u(x),Du(x),D^2u(x))\varphi^\alpha(x) dx=0
\end{align*}
be the weak formulation of a second order system of PDEs
\begin{align*}
f_\alpha(x,u(x),Du(x),D^2u(x))=0,\quad\alpha=1,2,...,m.
\end{align*}
Assume:
\begin{enumerate}
\item $K_f$ is invariant under a set of symmetry vector fields which span the tangent space $T_pE$ at each $p\in E$ (the invariance of $K_f$ is explained in \eqref{andin}).
\item Each symmetry vector field generates a corresponding conservation law of the form of a divergence expression (so-called continuity equation, see Definition \ref{contequ}).
\end{enumerate}
Then $f$ must be variational, that is, there exists a functional 
\begin{align*}
I(u)=\int_M L(x,u(x),Du(x),D^2u(x))dx
\end{align*}
with Lagrange function $L$, such that $K_f(u;\varphi)=\delta I(u;\varphi)$, where $\delta I$ is the first variation of $I$.\footnote{While in many classical cases the Lagrangian can be chosen to depend on $(x,u,Du)$ only, there are cases in which second derivatives $D^2u$ are needed, for example, the Monge-Ampère equation $u_{xx}u_{yy}-u_{xy}u_{xy}=0$.}
\end{theorem}
Later, the weak formulation $K_f$ will be replaced by the so-called source form $\Delta=f_\alpha du^\alpha\wedge dx$, see \eqref{30}, and we also consider the problem only locally, that is, for sufficiently small subsets $U\subset E$ and corresponding subset of $J^kE$, written as $(\pi^{k,0})^{-1}U$. We state our main theorem precisely in Section 4, but here we highlight in which ways it generalizes previously known results for second order system of PDEs:
\begin{itemize}
\item We do not make strong symmetry assumptions, like translation-, and gauge invariance in \cite{bib14,bib16}. We only assume that the symmetry vector fields span the tangent space $T_pE$ pointwise, see \eqref{48}.
\item As already mentioned, there is a connection between symmetries and conservation laws. Therefore, with the span-condition \eqref{48} we do not make strong conservation law assumptions, like divergence-free in \cite{bib14,bib16}. Notice that the divergence-free-condition forces a certain polynomial structure of the differential equations in the second order coordinates at the very beginning and, for example, the Monge-Ampère equation is excluded which is not the case in our theorem.
\item We do not need to assume that the set of symmetry vector fields, or certain subalgebras, are an abelian Lie algebra which is the case for translation invariance in \cite{bib14,bib16}. Our proof works for any non-abelian, or abelian Lie algebras which satisfy the mentioned span-condition \eqref{48}.
\item We do not need to assume that $n=m$, that is, the dimension of the base manifold $M$ is equal to the dimension of the fibers of $E$ which is assumed in \cite{bib14}. There are also other restrictions for the fibers of $E$ in \cite{bib16}. Our theorem works for any $n,m$ and without any restrictions on the fibers of $E$. For example, the vanishing of $H^{ab,i}_\Delta$ in \cite{bib14} is derived when explicitly using $n=m$.
\end{itemize}
In a sense, our result can be seen as the most general theorem for arbitrary system of PDEs of second order. Weakening our symmetry assumptions is only possible when restricting to first order systems which are free of topological obstructions, see \cite{bib15}, or having other restrictions, like equations which are polynomial in their arguments $u$ and derivatives $Du$, $D^2u$ and so on. Furthermore, we introduce a new inductive method to prove our main theorem. In our opinion, this method is much easier to handle and understand as, for example, the so-called $d$-fold method used in \cite[Lemma 2.3]{bib13} and \cite[p.12]{bib16} to solve similar problems. One reason for this is that the inductive method reduces a big problem to several very simple problems and the steps are easy to check. On the other hand, the $d$-fold method solves the mentioned problems in one relatively complicated step which is also hard to check. Of course, both methods can be useful in different situations and should be developed further.\\
\hspace*{0.5cm} Notice that the assumptions of symmetries and conservation laws in Theorem \ref{th0} are in some sense very natural, when we recall Noether\rq s theorem \cite{bib21}, since it can be seen as a reverse of Noether\rq s theorem. Roughly speaking, Noether\rq s theorems can be formulated as
\begin{align*}
&\begin{small}\text{'$f$ is variational' \quad \&\quad '$K_f$ is invariant' \quad $\Rightarrow$\quad '$f$ satisfies conservation laws'}\end{small},\\
&\begin{small}\text{'$f$ is variational' \quad \&\quad '$f$ satisfies conservation laws'\quad $\Rightarrow$\quad '$K_f$ is invariant'}\end{small}.
\end{align*}
Beside the both formulations of Noether\rq s theorems, we usually also distinguish the cases whether the symmetries are given by a finite dimensional group of transformations (also called Noether\rq s first theorems), or if they are given by an infinite dimensional group of transformations (also called Noether\rq s second theorems), but this does not make a major difference for the understanding of Takens\rq{} problem. Theorem \ref{th0}, that is, Takens\rq{} problem, is to prove 
\begin{align*}
\begin{small}\text{'$K_f$ is invariant' \quad \&\quad '$f$ satisfies conservation laws'\quad $\Rightarrow$\quad '$f$ is variational'}\end{small}.
\end{align*}
The paper is organized as follows. In Section 2 we present a simple instructive example which should help the reader to understand the problem in more detail, but without introducing much notation. In Section 3 we introduce our notation which is the notation of jet bundles, we precisely formulate what we mean by a differential equation, symmetry, conservation law, variational formulation, and we explain all the essential ideas behind Takens\rq{} problem. The main theorem and proof can be found in Section 4. Some of the longer calculations of Section 4 are carried out in Section 5 and 6. Finally, in Section 7 we provide a critical discussion of our result and explain the main open problems, especially in terms of applications in physics.\\
\quad\\
\textbf{Acknowledgment.} This article is a result of my Ph.D. thesis \cite{bib52}, written at the Carl von Ossietzky University Oldenburg, and I want to thank my supervisor Prof. Dr. Daniel Grieser very much for supporting this project with helpful advice and comments about this paper. In particular, I want to thank Daniel Grieser for giving me a lot of freedom in my research and supporting my interests. I also want to thank Prof. Dr. Ian M. Anderson for fruitful discussions about Takens\rq{} problem during my stay at Utah State University and for a very pleasant time.

\section{A simple instructive example}
\label{sec:1}
In this example we consider the differential equation $f=0$ itself, instead of the weak formulation $\int f\varphi dx=0$. This is because we then need to introduce less notation, and for translations, $f$ and $K_f$ satisfy the same symmetry and both formulations are in some sense equivalent. Let $f(x,u(x),u_x(x),u_{xx}(x))=0$ be a scalar differential equation for the unknown function $u:\mathbb{R}\to\mathbb{R}$, where $f=u+u_{xx}$ describes the harmonic oscillator. As usual, $u_x(x)=\frac{du(x)}{dx}$ and $u_{xx}(x)=\frac{d^2u(x)}{dx^2}$. This equation is translation invariant in $x$-direction, that is, $f(x+s,u,u_x,u_{xx})=f(x,u,u_x,u_{xx})$ for all $s\in\mathbb{R}$ and for all $(x,u,u_x,u_{xx})$. The infinitesimal condition for this symmetry is given by the equation
\begin{align*}
\frac{d}{ds}f(x+s,u,u_x,u_{xx})\vert_{s=0}=\frac{\partial}{\partial x}f(x,u,u_x,u_{xx})=0\quad\text{for all $(x,u,u_x,u_{xx})$},
\end{align*}
with symmetry vector field $\frac{\partial}{\partial x}$. When we consider $x$ as the time variable, and $u$ as the (one-dimensional) position of a particle at the time $x$, then the equation $f=u+u_{xx}=0$ satisfies energy conservation. More precisely, when we multiply $f$ by $u_x$, then we get
\begin{align}
u_xf=u_x(u+u_{xx})=\frac{d}{dx}(\frac{1}{2}u^2_x+\frac{1}{2}u^2)=0\quad&\text{on solutions, where $f=0$}\nonumber\\
\Rightarrow\quad\frac{1}{2}(u^2_x+u^2)=\text{constant}\quad&\text{on solutions, where $f=0$.}\label{37}
\end{align}
The term $\frac{1}{2}(u^2_x+u^2)$ can be interpreted as the kinetic plus the potential energy of the particle and it delivers a conservation law, called energy conservation. As we already mentioned above, the structure of the theorem we prove can be written as
\begin{align*}
\begin{small}\text{'$f$ (or $K_f$) is invariant' \quad \&\quad '$f$ satisfies conservation laws'\quad $\Rightarrow$\quad '$f$ is variational'}\end{small}.
\end{align*} 
The fact that $f=u+u_{xx}$ is variational is easy to check in this case. For example, let us consider the Lagrangian $L=\frac{1}{2}u^2-\frac{1}{2}u^2_x$, and we get $\frac{\partial L}{\partial u}-\frac{d}{dx}\frac{\partial L}{\partial u_x}=u+u_{xx}=f$. It is well-known, from Noether\rq s first theorem, that time translation invariance leads to energy conservation, and this correlation can be observed in our example here. Actually, so far, in this example we did not have to prove anything, but now we may state the following conjecture.\\
\quad\\
\textit{Any scalar differential equation $f(x,u(x),u_x(x),u_{xx}(x))=0$ of second order, where $u:\mathbb{R}\to\mathbb{R}$, which is $\frac{\partial}{\partial x}$-invariant, and which has a corresponding conservation law of the form $u_xf=\frac{d}{dx}E$, where $E$ is a no more specified (energy) function, can be written as the Euler-Lagrange equation for some Lagrangian $L$. (Let us also assume that the functions $f$ and $E$ are smooth in $x,u,u_x$ and so on.)}\\
\quad\\
The conjecture is actually true, as one can check in the Appendix, and it does not require a lot of non-standard notation to prove it in this scalar ODE case. However, the proof in the Appendix does not provide a systematic algorithm how to solve such problems in general, and the proof of our main theorem will be different. Hopefully, the proof in the Appendix helps the reader to get a better understanding how complex it can be to prove such theorems systematically.

\section{Preliminaries}
\label{sec:2}
Let $x=(x^i)^{i=1,2,...,n}$ denote the independent coordinates, and $u=(u^\alpha)^{\alpha=1,2,...,m}$ the dependent coordinates which are used to describe a $k$-th order system of PDEs $f(x,u(x),Du(x),...,D^ku(x))=0$, where we use the short notation
\begin{align*}
D^lu(x):=\left(\frac{\partial^l u^\alpha(x)}{\partial x^{i_1}\partial x^{i_2}...\partial x^{i_l}}\right)^{\alpha=1,2,...,m}_{i_1,i_2,...,i_l=1,2,...,n},\quad 1\le l\le k,
\end{align*}
for the set of all possible $l$-th order partial derivatives of $u(x)$. When we want to consider PDEs on manifolds, or when we want to make our descriptions independent of the choice of coordinates, then it is reasonable to use the language of fiber- and jet bundles.

\paragraph{The jet bundle.} Let $\pi:E\to M$ be a fiber bundle with local coordinates $x$ on $M$ and adapted local coordinates $(x,u)$ on $E$, that is, $u$ describes the coordinates of the fibers of $E$. We have $\operatorname{dim}E=n+m$ and $\operatorname{dim}M=n$. Local charts on $E$ are written as $(U,\varphi)$, where $U\subset E$ is open, and we have $\varphi:U\to\mathbb{R}^{n+m}$. Furthermore we write $(\bar U,\varphi^0)$ for local charts on $M$, where $\bar U\subset M$ is open, and we have $\varphi^0:\bar U\to\mathbb{R}^n$. Adapted local coordinates means that $\bar U=\pi U$ and that $\varphi,\varphi^0$ satisfy $\varphi^0(\pi p)=\tilde \pi\varphi(p)$ for all $p\in U$, where $\tilde\pi:\mathbb{R}^{n+m}\to\mathbb{R}^n$ is the canonical projection.\\
\hspace*{0.5cm} (Local) sections of $\pi$ are written as $\sigma:\Omega\to E$, where $\Omega\subset M$, and $\pi\sigma$ is the identity map on $\Omega$. We have $\varphi\circ \sigma\circ (\varphi^0)^{-1}(x)=(x,s(x))$ on $\varphi^0(\Omega\cap\bar U)$, and in local coordinates the section $\sigma$ can be identified with the function $s:\mathbb{R}^n\to \mathbb{R}^m$. In fiber bundle notation, a differential equation is written as $f(x,s(x),Ds(x),...,D^ks(x))=0$, where $f$ is a function defined on the so-called jet bundle, and solutions are certain sections of $\pi$.\\
\hspace*{0.5cm} The $k$-th order jet bundle of $E$ is written as $J^kE$ and it has (adapted) local coordinates 
\begin{align*}
(x^i,u^{\alpha_0},u^{\alpha_1}_{i_1},...,u^{\alpha_k}_{j_1j_2...j_k})^{\alpha_0,\alpha_1,...,\alpha_k=1,2,...,m}_{i,i_1,...,j_1,...,j_k=1,2,...,n},
\end{align*}
where the indices $i_1i_2...i_l$ in the coordinates $u^\alpha_{i_1i_2...i_l}$ are ordered as $1\le i_1\le i_2\le...\le i_l\le n$ for every $1\le l\le k$. These local coordinates are in one to one correspondence with a $k$-th order equivalence class $[\sigma]_k(q)$ of local sections of $\pi$ at a point $q\in M$, where $\varphi^0(q)=x$, such that all local sections in this equivalence class satisfy (by the definition of equivalence class)
\begin{align}
u^\alpha&=s^\alpha(\tilde x)\vert_{\tilde x=x},\quad\quad\quad\quad\alpha=1,2,...,m,\nonumber\\
u^\alpha_{i_1i_2...i_l}&=\frac{\partial^l s^\alpha(\tilde x)}{\partial \tilde x^{i_1}...\partial \tilde x^{i_l}}\vert_{\tilde x=x},\quad\;\;\,\alpha=1,2,...,m,\quad 1\le i_1\le i_2\le ...\le i_l\le n,\label{18}
\end{align}
for every $1\le l\le k$. The equivalence classes $[\sigma]_k(q)$ are the points in $J^kE$. For further details see \cite{bib9} and \cite{bib36}. The ordering $1\le i_1\le i_2\le ...\le i_l\le n$ is reasonable when defining the local coordinates $u^\alpha_{i_1i_2...i_l}$, for example, to get the correct dimension of $J^kE$, but from now on we also allow expressions $u^\alpha_{i_1i_2...i_l}$ to occur, where $i_1i_2...i_l$ are not ordered in an increasing manner which makes the calculations below much easier. Whenever such an expression occurs, we make the identification $u^\alpha_{i_1...i_l}=u^\alpha_{\tau(i_1)...\tau(i_l)}$ for any permutations $\tau$ of $i_1...i_l$. For example, $u^\alpha_{12}=u^\alpha_{21}$.\footnote{Later we will introduce the so-called total derivative operator $D_i$ and it satisfies $D_iu^\alpha_{i_1i_2...i_l}=u^\alpha_{ii_1i_2...i_l}$ when we use the identification and we do not have to care in which position we have to write the indice $i$. Another reason is that the identification also comes out naturally by \eqref{18}, since partial derivatives commute.} Furthermore, we use Einstein summation whenever there occur repeated lower and upper indices. Latin indices run from $1$ to $n$ and greek indicies from $1$ to $m$. This means in particular that each of the indices $i_1,...,i_l$ in an expression $u^\alpha_{i_1...i_l}A^{i_1...i_l}_\alpha$ runs from $1$ to $n$ and there is no ordering of $i_1...i_l$ in the Einstein summation.\\ 
\hspace*{0.5cm} Let us introduce the short notation $D^lu:=(u^\alpha_{i_1i_2...i_l})^{\alpha=1,2,...,m}_{i_1,i_2...,i_l=1,2,...n}$, $1\le l\le k$ for higher order coordinates on $J^kE$. Then a point in $J^kE$ is given by the coordinates $(x,u,Du,...,D^ku)$. We call $(x,u,Du,...,D^k u)$ the \textbf{$k$-th order jet coordinates} and a single coordinate $u^\alpha_{i_1i_2...i_l}$, $0\le l\le k$ is called \textbf{$l$-th order coordinate} (and $(x,u)$ are the 0-th order jet coordinates).\\
\hspace*{0.5cm} The jet bundle has different types of projections:
\begin{align*}
\pi^{k}&:J^kE\to M,\\
\pi^{k,0}&:J^kE\to E,\\
\pi^{k,l}&:J^kE\to J^lE,\quad 0\le l\le k,
\end{align*}
where we define $J^0E=E$ and $\pi^{k,k}$ is the identity map on $J^kE$. These projections are an immediate consequence of the definition of $J^kE$, since we can always map a point $[\sigma]_k(q)\in J^kE$ to a point  $\pi^{k,l}[\sigma]_k(q)=[\sigma]_l(q)\in J^lE$, $l\le k$, when ignoring the equivalence of higher order derivatives of sections $s$ in the equivalence class $[\sigma]_k(q)$, see \eqref{18}.\\
\hspace*{0.5cm} If $g\in C^\infty(J^kE)$ then we call the number
\begin{align}
l=\min\limits_{\substack{h\in C^\infty(J^rE):\\g=(\pi^{k,r})^\ast h}}\{r\},\label{order}
\end{align}
the \textbf{order of $g$}. For example, $u^\beta_iu^\gamma_{kl}$ is a second order function.

\paragraph{Prolongation of sections.} As already mentioned, points in $J^kE$ are given as equivalence classes of sections of $\pi$. There is a very special set of points in $J^kE$, namely the set which consists of a so-called prolonged section. Let $\sigma:\Omega\to E$, $\Omega\subset M$ be a (local) section of $\pi$. Then we define $\text{j}^k\sigma:\Omega\to J^kE$ as the (local) section of $\pi^k$, which maps $q\mapsto [\sigma]_k(q)$ for every $q\in\Omega$. This map is called \textbf{jet prolongation of $\sigma$}, see \cite[p.30]{bib9}. Prolongation lifts a section of $\pi$ to a section of $\pi^k$. Notice that every point $[\sigma]_k(q)\in J^kE$ can also be written as $\text{j}^k\sigma(q)$ for a certain section $\sigma$ of $\pi$ and a certain $q\in M$. Sometimes it is appropriate to write $\text{j}^k\sigma(q)$ and sometimes $[\sigma]_k(q)$.

\paragraph{Vector fields.} Let $\mathfrak{X}(J^kE)$ denote the set of vector fields on $J^kE$. In local coordinates, a vector field $V\in\mathfrak{X}(J^kE)$ is written as\footnote{When using concrete labels, like $V^1$, then we write $V^{x^1}$, or $V^{u^1}$ to indicate if it belongs to a coefficient $V^i$, or $V^\alpha$.}
\begin{align*}
V=V^i\frac{\partial}{\partial x^i}+V^\alpha\frac{\partial}{\partial u^\alpha}+V^\alpha_{i_1}\frac{\partial}{\partial u^\alpha_{i_1}}+...+V^\alpha_{i_1i_2...i_k}\frac{\partial}{\partial u^\alpha_{i_1i_2...i_k}}.
\end{align*}
Notice that we may assume that $V^\alpha_{i_1...i_l}=V^\alpha_{\tau(i_1)...\tau(i_l)}$ for every permutation $\tau$, since by our previous definition we have $\frac{\partial}{\partial u^\alpha_{i_1i_2...i_l}}=
\frac{\partial}{\partial u^\alpha_{\tau(i_1)...\tau(i_l)}}$, and we may assume that the coefficients $V^\alpha_{i_1i_2...i_l}$ are completely symmetrized in the indices $i_1i_2...i_l$. There are $\frac{l!}{l_1!l_2!...l_n!}$ possibilities how to rearrange the numbers $i_1,i_2,...,i_l$, where $l_r$ is the number of occurrences of the index $r$ in the multi-index $i_1i_2...i_l$ and we have $1\le r\le n$. Therefore, let us also define the weighted partial derivatives
\begin{align*}
\partial_i:=\frac{\partial}{\partial x^i},\quad\partial_\alpha:=\frac{\partial}{\partial u^\alpha},\quad \partial^{i_1i_2...i_l}_\alpha:=\frac{l_1!l_2!...l_n!}{l!}\frac{\partial}{\partial u^\alpha_{i_1i_2...i_l}},\quad 1\le l\le k,
\end{align*}
see \cite{bib15}.\footnote{Again, when we have concrete labels, like $\partial_1$, then we write $\partial_{x^1}$ or $\partial_{u^1}$ to indicate if it belongs to $\partial_i$ or $\partial_\alpha$.}\\
\hspace*{0.5cm} On fiber bundles, there are certain vector fields with additional structure, given by the projection map $\pi$. We call a vector field $V\in\mathfrak{X}(E)$ \textbf{$\pi$-projectable}, if $\pi_\ast V$ exists, and if it is a vector field on $M$. Moreover, we call a vector field $V\in\mathfrak{X}(E)$ \textbf{$\pi$-vertical}, if $\pi_\ast V=0$.
In local coordinates, these vector fields are written as
\begin{align*}
&V=V^i(x)\partial_i+V^\alpha(x,u)\partial_\alpha,\,\;\quad\text{($\pi$-projectable vector field on $E$)},\\
&V=V^\alpha(x,u)\partial_\alpha,\quad\quad\quad\quad\quad\quad\text{($\pi$-vertical vector field on $E$)}.
\end{align*}
Notice that the flow of a $\pi$-projectable vector field transforms sections of $\pi$ to sections of $\pi$ and for this reason they will be very important later on. Since $J^kE$ has different types of projections $\pi^k$, $\pi^{k,l}$, we also define $\pi^k$-, \textbf{$\pi^{k,l}$-projectable} and $\pi^k$-, \textbf{$\pi^{k,l}$-vertical} vector fields in the same way.\\
\hspace*{0.5cm} There is another kind of vector field we want to mention here and this is the so-called total derivative. The \textbf{total derivative on $J^kE$}, written with respect to some local coordinates on $J^{k+1}E$, is defined as
\begin{align}
D_i=D_i(k):=\partial_i+u^\alpha_i\partial_\alpha+u^\alpha_{ii_1}\partial^{i_1}_\alpha
+...+u^\alpha_{ii_1i_2...i_k}\partial^{i_1i_2...i_k}_\alpha,\label{divop}
\end{align}
and it is a kind of vector field, however it is not vector field in the usual sense, since the last coefficients $u^\alpha_{ii_1i_2...i_k}$ are not defined in $J^kE$, rather in $J^{k+1}E$. A coordinate invariant definition is the following one. Let $V=V^i\partial_i$ be a vector field on $M$ and $\phi^0_t$ its flow on $M$. Furthermore, let $[\sigma]_{k+1}(q)$ be a point in $J^{k+1}E$. Then we assign a tangent vector $\text{tot}^kV$ at the point $\pi^{k+1,k}[\sigma]_{k+1}(q)=[\sigma]_k(q)\in J^kE$ as follows. We take any local section $\sigma:\Omega\to E$ in the equivalence class $[\sigma]_{k+1}(q)$ and we compute the tangent vector
\begin{align*}
\text{tot}^kV\vert_{[\sigma]_k(q)}:=\frac{d}{dt}\text{j}^k\sigma(\phi^0_t(q))\vert_{t=0}
\end{align*}
at the point $[\sigma]_k(q)$ which always works for sufficiently small $t$, such that we have $\phi^0_t(q)\in \Omega$. The tangent vector $\text{tot}^kV\vert_{[\sigma]_k(q)}$ depends on the equivalence class $[\sigma]_{k+1}(q)$. When we consider the union of all these tangent vectors at all points $[\sigma]_k(q)\in J^kE$ then we get a kind of vector field on $J^kE$. However, this is not a vector field in the usual sense as a map $J^kE\to TJ^kE$, it is rather a map $J^{k+1}E\to TJ^kE$ (sometimes called a vector field along the map $\pi^{k+1,k}$, or a generalized vector field). In local coordinates, we get $\text{tot}^kV=V^iD_i$, where we again find the total derivative on $J^kE$. We call $\text{tot}^kV$ the \textbf{total vector field of $V$ on $J^kE$}. More generally, for every $\pi^k$-projectable vector field $W$ on $J^kE$ we can assign a unique total vector field $\text{tot}^kW$ by the above construction, where  $V=\pi^{k}W$ is the vector field on $M$. Sometimes this construction is also called horizontalization of $W$, see \cite[p.17]{bib36}. For later purposes it is important to notice here that when we apply a total vector field to functions we get a map $\text{tot}^kV:C^\infty(J^kE)\to C^\infty (J^{k+1}E)$, that is, the order is increased by one. This also means that when we write several total derivatives repeated, like $D_iD_j$, we actually mean $D_i(k+1)D_j(k)$. Fur further details also see \cite{bib7,bib20}.

\paragraph{Prolongation of flows and vector fields.} As we saw above, a section of $\pi$ can be lifted to a section of $\pi^k$ and we called this lift jet prolongation of the section. In a similar way we can lift a $\pi$-projectable vector field on $E$ to a vector field on $J^kE$. Let $V$ be a $\pi$-projectable vector field on $E$ and $\phi_t$ its flow. Moreover, let $\pi_\ast V$ be the corresponding vector field on $M$ with flow $\phi^0_t$. It can be shown that $\phi_t\circ\sigma\circ(\phi^0_t)^{-1}$ is a section of $\pi$ (basically by showing that $\pi(\phi_t\circ\sigma\circ(\phi^0_t)^{-1})(q)=q$ for all $q\in\Omega\subset M$). Now we can prolong this section and we get a section of $\pi^k$. Let us define $\text{j}^k\phi_t$ as the \textbf{jet prolongation of the flow $\phi_t$} through
\begin{align*}
\text{j}^k\phi_t(\text{j}^k\sigma(q)):=\text{j}^k[\phi_t\circ\sigma\circ(\phi^0_t)^{-1}](\phi^0_t(q))
\end{align*}
for every point $\text{j}^k\sigma(q)\in J^kE$. We get a map $\text{j}^k\phi_t:J^kE\to J^kE$ and it is also an isomorphism, see \cite[p.32]{bib9}. The section $\phi_t\circ\sigma\circ(\phi^0_t)^{-1}$ depends on the parameter $t$, and when we prolong it and take the derivative with respect to $t$
\begin{align*}
\text{j}^kV(\text{j}^k\sigma(q)):=\frac{d}{dt}\text{j}^k[\phi_t\circ\sigma\circ(\phi^0_t)^{-1}](\phi^0_t(q))\vert_{t=0},
\end{align*}
we get a vector field $\text{j}^kV$ on $J^kE$ at $t=0$, depending on the point $\text{j}^k\sigma(q)\in J^kE$ which is called the \textbf{jet prolongation of the vector field $V$}. Such vector fields are used to describe symmetries of differential equations, or to describe perturbations of integral functionals in the calculus of variations. In local coordinates, a prolonged vector field can be written as
\begin{align}
\text{j}^kV=&V^i\partial_i+V^\alpha \partial_\alpha+\xi^\alpha_{i_1}\partial^{i_1}_\alpha+\xi^\alpha_{i_1i_2}\partial^{i_1i_2}_\alpha+...+\xi^\alpha_{i_1i_2...i_k}\partial^{i_1i_2...i_k}_\alpha,\label{n58}
\end{align}
where the $\xi$-coefficients are defined recursively as 
\begin{align*}
\xi^\alpha_{i_1}&=D_{i_1}V^\alpha-u^\alpha_iD_{i_1}V^i,\\
\xi^\alpha_{i_1i_2...i_l}&=D_{i_l}\xi^\alpha_{i_1i_2...i_{l-1}}-u^\alpha_{ii_1i_2...i_{l-1}}D_{i_l}V^i,\quad 2\le l\le k,
\end{align*}
see \cite[p.32]{bib9} and \cite[p.26]{bib36}. It also turns out, that on $TJ^{k+1}E$, the prolongation $\text{j}^kV$ can be written equivalently as
\begin{align}
\text{j}^kV=
V^iD_i+V^\alpha_{ch} \partial_\alpha+(D_{i_1}V^\alpha_{ch})\partial^{i_1}_\alpha+...+
(D_{i_1}D_{i_2}...D_{i_k}V^\alpha_{ch})\partial^{i_1i_2...i_k}_\alpha,\label{n49}
\end{align}
where $(V^\alpha_{ch}):=(V^\alpha-u^\alpha_iV^i)$ is called \textbf{characteristic of $V$}, and $V_{ch}=V^\alpha_{ch}\partial_\alpha$ is the \textbf{characteristic vector field of $V$}, see \cite[p.118]{bib0}. The right hand side in \eqref{n49} is not defined in $TJ^kE$, rather as a map $J^{k+1}E\to TJ^kE$, the same as with the total derivative operator. If there is no danger of confusion we either write \eqref{n58}, or \eqref{n49} for $\text{j}^kV$ and we do not explicitly notice one or the other. The formal prolongation of $V_{ch}$ is defined as
\begin{align*}
\text{j}^kV_{ch}:=V^\alpha_{ch}\partial_\alpha+(D_iV^\alpha_{ch})\partial^i_\alpha+...+
(D_{i_1}D_{i_2}...D_{i_k}V^\alpha_{ch})\partial^{i_1i_2...i_k}_\alpha.
\end{align*}
Then \eqref{n49} can be written as $\text{j}^kV=\text{tot}^kV+\text{j}^kV_{ch}$. Roughly speaking, $\text{j}^kV$ decomposes into a total-, and a $\pi^k$-vertical part, but the decomposition cannot be done in $TJ^kE$, it rather has to be considered as a map $J^{k+1}E\to TJ^kE$.\\
\hspace*{0.5cm} When we think of $\pi$-projectable vector fields (or their prolongations) as symmetries and Lie algebras, then the following property is very useful. If $V,W$ are $\pi$-projectable vector fields on $E$ and if they form a Lie algebra, then $\text{j}^kV,\text{j}^kW$ also form a Lie algebra and we have $[\text{j}^kV,\text{j}^kW]=\text{j}^k[V,W]$, where $[.,.]$ is the Lie bracket, see \cite[p.29]{bib36}.

\paragraph{Zero-functionals and the Lagrange form.} In the calculus of variations we investigate integral functionals of the form
\begin{align}
I_{\sigma}=\int_{\Omega} (\text{j}^k\sigma)^\ast (Ldx)=\int_{\Omega} L(x,s(x),Ds(x),...,D^ks(x))dx,\label{25}
\end{align}
where $L=L(x,u,...,D^{k}u)$ is the Lagrangian, the closure $\bar \Omega\subset M$ is a compact set, and $\Omega$ is oriented. We use the short notation $dx:=dx^1\wedge dx^2\wedge ...\wedge dx^n$. We also call \eqref{25} a \textbf{zero-functional}. The so-called \textbf{Lagrange form}
\begin{align*}
\lambda:=Ldx
\end{align*}
is a main ingredient in the integral functional \eqref{25}, and for many calculations it is sufficient to consider $\lambda$ instead of \eqref{25}. For an intrinsic definition of $\lambda$ and further information see \cite{bib36} in $\S 4$.\\
\hspace*{0.5cm} Let $\phi_t$ be the flow of a $\pi$-projectable vector field $V\in\mathfrak{X}(E)$ and $\phi^0_t$ the flow of $\pi_\ast V$. The \textbf{transformation of zero-functionals under $\phi_t$} is defined as (a formal pullback)
\begin{align}
\phi^\ast_t I_{\sigma}:=\int_{\phi^0_t\Omega}\text{j}^k(\phi_t\circ\sigma\circ\phi^0_{-t})^\ast\lambda=\int_\Omega(\text{j}^k\sigma)^\ast[(\text{j}^k\phi_t)^\ast\lambda],\label{36}
\end{align}
see \cite[p.42]{bib9} and \cite[p.111]{bib36} (we may assume that the transformation is defined for at least sufficiently small $t$). In the following, $\mathcal{L}_W$ denotes the Lie derivative with respect to the vector field $W$. The expression \eqref{36} allows us to define a \textbf{formal Lie derivative of zero-functionals} as
\begin{align}
\mathcal{L}_{\text{j}^kV}I_{\sigma}:=\frac{d}{dt}\int_{\phi^0_t\Omega}\text{j}^k(\phi_t\circ\sigma\circ\phi^0_{-t})^\ast\lambda\vert_{t=0}=\int_{\Omega}(\text{j}^k\sigma)^\ast(\mathcal{L}_{\text{j}^kV}\lambda),\label{27}
\end{align}
where $\mathcal{L}_{\text{j}^kV}\lambda$ is the Lie derivative in the usual sense, applied to the differential form $\lambda$, see \cite[p.42]{bib9} and \cite[p.111]{bib36}. The formal Lie derivative in \eqref{27} can be used for two things: 
\begin{enumerate}
\item When $V\in\mathfrak{X}(E)$ is a $\pi$-vertical vector field, such that $\operatorname{supp}V\subset\pi^{-1}\Omega$, then \eqref{27} describes the first variation of $I$ in direction of $V$. Support in $\pi^{-1}\Omega\subset E$ implies that $V$ must vanish at the boundary $\pi^{-1}\partial\Omega\subset E$ and this means that we can do partial integration without getting boundary terms which will lead to the Euler-Lagrange equation.
\item Symmetries of $I$ are defined as $\pi$-projectable vector fields $V\in\mathfrak{X}(E)$, such that \eqref{27} vanishes for all sections $\sigma$ of $\pi$.
\end{enumerate}

\paragraph{One-functionals and the source form.} Let $R\in\mathbb{N}$ and
\begin{align}
(\text{j}^k\sigma)^\ast f_r=f_r(x,s(x),Ds(x),...,D^ks(x))=0,\quad r=1,2,...,R\label{28}
\end{align}
be a system of PDEs, where $(f_r)$ are functions defined on $J^kE$, and the section $\sigma$ is a solution of the differential equation. We can multiply \eqref{28} by test functions $\varphi_r=\varphi_r(x)$ which are defined on $\Omega\subset M$, the closure $\bar\Omega$ is compact, and we can integrate over $\Omega$, such that we get
\begin{align}
\tilde K_{\sigma}(\varphi)=\int_\Omega \sum_{r=1}^{R}(\text{j}^k\sigma)^\ast(f_r\varphi_rdx)=0,\quad\text{(for all test functions)}
\end{align}
which is known as a \textbf{weak formulation of \eqref{28}}. As usual, test functions are $C^\infty(M)$-functions with support in $\Omega$. When $R=m$, then we can define a very special weak formulation, namely
\begin{align}
\tilde K_{\sigma}(V)=\int_\Omega (\text{j}^k\sigma)^\ast(f_\alpha V^\alpha dx)=0,\label{29}
\end{align}
where the test functions $\varphi_r$ can be identified with $\pi$-vertical vector fields $V=V^\alpha\partial_\alpha\in\mathfrak{X}(E)$, and they have the additional property that their support is in $\pi^{-1}\Omega\subset E$, such that $(\text{j}^k\sigma)^\ast V^\alpha$ is a test function with support $\Omega$. The Einstein summation in \eqref{29} shows that $(f_\alpha)$ now must have a very specific transformation property under local coordinate changes which is not necessarily the case for any equation of the form \eqref{28}. For example, let $L=L(x,u,Du)$ be a first order Lagrangian, then the first variation \eqref{27} leads to such a particular weak formulation, where $f_\alpha=(\partial_\alpha-D_i\partial^i_\alpha)L$ is the so-called Euler-Lagrange expression which satisfies this transformation property. We can also work with \eqref{29} when $(f_\alpha)$ cannot be written as an Euler-Lagrange expression, but the transformation property for $(f_\alpha)$ must be the same as for Euler-Lagrange equations. More precisely, it is given as $f_\alpha\to f_\alpha\frac{\partial u^\alpha}{\partial v^\beta}\det\frac{\partial x}{\partial y}$, where $\frac{\partial x}{\partial y}$ is a short notation for the Jacobian matrix and the $\det\frac{\partial x}{\partial y}$-term comes from the differential $dx$ in \eqref{29}.\\
\hspace*{0.5cm} As we have seen above, that the main ingredient of zero-functionals is the Lagrange form $\lambda$, we now want to find, in a similar way, the main ingredient of the weak formulation \eqref{29}. It turns out that this is the so-called \textbf{source form}
\begin{align}
\Delta:=f_\alpha du^\alpha\wedge dx,\label{30}
\end{align}
which satisfies the above mentioned transformation property. The source form is a $(n+1)$-form on $J^kE$. For an intrinsic definition and further details see \cite{bib17} and \cite[p.560]{bib45}.\footnote{Because $R=m$, and because of the specific transformation of $(f_\alpha)$, we cannot assign a source form to every differential equation of the form \eqref{28} and the assignment is also not unique. For example, let $R=2$ and assume $(f_r)$ is given,  then we can assign $\Delta=f_1du^1\wedge dx+f_2du^2\wedge dx$, or $\Delta=f_2du^1\wedge dx+f_1du^2\wedge dx$. Also see the discussion in Section 7.} With the help of the source form, we define the \textbf{one-functional}
\begin{align}
K_{\sigma}(W):=\int_\Omega (\text{j}^k\sigma)^\ast(W\lrcorner\Delta),\label{32}
\end{align}
for every $W\in\mathfrak{X}(J^kE)$. Then \eqref{29} can be written equivalently as $K_\sigma(\text{j}^kV)=0$ for all $\pi$-vertical $V\in\mathfrak{X}(E)$ with $\operatorname{supp}V\subset\pi^{-1}\Omega$.\\ 
\hspace*{0.5cm} Let $\phi_t$ be the flow of a $\pi$-projectable vector field $V\in\mathfrak{X}(E)$. As before, we define the \textbf{transformation of one-functionals under $\phi_t$} as (a formal pullback)
\begin{align}
(\phi^\ast_tK)_\sigma(W)&:=\int_{\Omega}(\text{j}^k\sigma)^\ast\{W\lrcorner[(\text{j}^k\phi_t)^\ast\Delta]\}\label{dist}
\end{align}
for every section $\sigma$ and for every $W\in\mathfrak{X}(J^kE)$. Notice that \eqref{dist} has to be distinguished from
\begin{align}
\phi^\ast_t(K_\sigma(W))&:=\int_{\Omega}(\text{j}^k\sigma)^\ast[(\text{j}^k\phi_t)^\ast(W\lrcorner\Delta)]\label{dist2}
\end{align}
which is the transformation of a certain zero-functional. Moreover, we define the \textbf{formal Lie derivative of one-functionals} as
\begin{align}
(\mathcal{L}_{\text{j}^kV}K)_\sigma(W)&:=\frac{d}{dt}(\phi^\ast_tK)_\sigma(W)\vert_{t=0}=\int_\Omega(\text{j}^k\sigma)^\ast[W\lrcorner(\mathcal{L}_{\text{j}^kV}\Delta)]\label{andin}
\end{align}
for every section $\sigma$ and for every $W\in\mathfrak{X}(J^kE)$, where $V\in\mathfrak{X}(E)$ is $\pi$-projectable and $\phi_t$ the flow of $V$. On the other hand, \eqref{dist2} leads to
\begin{align}
\mathcal{L}_{\text{j}^kV}(K_\sigma(W))&:=\frac{d}{dt}\phi^\ast_t(K_\sigma(W))\vert_{t=0}=\int_\Omega(\text{j}^k\sigma)^\ast[\mathcal{L}_{\text{j}^kV}(W\lrcorner\Delta)].\label{andin2}
\end{align}
Formula \eqref{andin} will be used when we define symmetries, and \eqref{andin2} is helpful when we derive the so-called integrability conditions, what we will do next. However, both formulas are connected through the identity $(\mathcal{L}_{\text{j}^kV} K)(W)-\mathcal{L}_{\text{j}^kV}(K(W))=K([W,\text{j}^kV])$.

\paragraph{Integrability conditions.} We say that a one-functional $K$ is \textbf{variational}, if there exists a zero-functional $I$, such that 
\begin{align}
\mathcal{L}_{\text{j}^kV}I_\sigma=K_\sigma(\text{j}^kV)\label{notice}
\end{align}
is satisfied for all sections $\sigma$, and for all $\pi$-vertical $V\in\mathfrak{X}(E)$ with $\operatorname{supp}V\subset\pi^{-1}\Omega$. This is analogous how we define an exact one-form $\omega$, where a zero-form $\eta$ must exist, such that the equation 
\begin{align}
(\mathcal{L}_v\eta)_p=(v\lrcorner\omega)_p\label{noticeb}
\end{align}
is satisfied for all points $p$, and for all vector fields $v$. The only difference is that \eqref{notice} only holds for certain vector fields with the additional assumptions \textbf{'$\pi$-vertical'} and \textbf{'support in $\pi^{-1}\Omega\subset E$'}, whereas \eqref{noticeb} holds for \textbf{all} vector fields, and that the points in the calculus of variations are chosen from an \textbf{infinite dimensional space} of sections, whereas the points $p$ in \eqref{noticeb} are chosen from a \textbf{finite dimensional space}.\\
\hspace*{0.5cm} By partial integration, the definition in \eqref{notice} can be transferred to source forms and to the coefficients $(f_\alpha)$, where we get the following definition.
\begin{definition}\label{definitionVar}
Let $\Delta=f_\alpha du^\alpha\wedge dx$ be a source form on $J^kE$. We say that $\Delta$ is \textbf{locally variational}, if for each $p\in E$ there exists a neighborhood $U\subset E$, a Lagrange form $\lambda =Ldx$ on $(\pi^{k,0})^{-1}U\subset J^kE$, and local coordinates on $(\pi^{k,0})^{-1}U\subset J^kE$, such that we can write
\begin{align}
f_\alpha=(\partial_\alpha-D_{i_1}\partial^{i_1}_\alpha+...
+(-1)^kD_{i_1}D_{i_2}...D_{i_k}\partial^{i_1i_2...i_k}_\alpha)L\label{31}
\end{align}
on $(\pi^{k,0})^{-1}U$. If there exists a Lagrange form $\lambda$ on $J^kE$, such that we can write $\Delta$ as \eqref{31} in every local coordinate system, then we say $\Delta$ is \textbf{globally variational}.
\end{definition}
Notice that \eqref{31} is independent of the choice of local coordinates. For a coordinate independent definition also see \cite{bib32}.\\
\hspace*{0.5cm} Now we want to derive conditions for $(f_\alpha)$ under which they can be written as \eqref{31} for some $L$. We start the discussion with a simple case. Let us assume that we have a source form $\Delta=f_\alpha du^\alpha\wedge dx$ of $0$-th order, that is, where $(f_\alpha)$ only depend on $(x,u)$. Then it is easy to check that $\Delta$ is locally variational if there exists a Lagrange form $\lambda=Ldx$ on $E$, such that locally $\Delta=d\lambda$. A necessary and sufficient condition for a 0-th order source form $\Delta$ to be locally variational is that $d\Delta=0$, that is, $\Delta$ must be closed and this leads to the integrability conditions $\partial_\beta f_\alpha(x,u)-\partial_\alpha f_\beta(x,u)=0$. This is well-known by the theory of locally exact sequences. For higher order source forms, that is, where $(f_\alpha)$ depend on $(x,u,Du,...,D^k u)$, the necessary and sufficient conditions for $\Delta$ to be locally variational are the so-called \textbf{Helmholtz conditions}. When $\Delta$ is of second order, then they are given by the equations
\begin{align}
&\partial_\beta f_\alpha-\partial_\alpha f_\beta+D_i\partial^i_\alpha f_\beta-D_iD_j\partial^{ij}_\alpha f_\beta
=0,\nonumber\\
&\;\;\quad\quad\partial^i_\beta f_\alpha+\partial^i_\alpha f_\beta-2D_j\partial^{ij}_\alpha f_\beta
=0,\nonumber\\
&\;\;\quad\quad\quad\quad\quad\partial^{ij}_\beta f_\alpha-\partial^{ij}_\alpha f_\beta
=0.\label{h}
\end{align}
These conditions are known for a very long time and there is quite a lot of literature about them \cite{bib40,bib9,bib36,bib29}. Beside the Helmholtz conditions, we also define the so-called \textbf{Helmholtz expressions} (for second order $\Delta$) as
\begin{align}
\left(\begin{matrix}H_{\alpha\beta}(\Delta)\\ H^i_{\alpha\beta}(\Delta)\\ H^{ij}_{\alpha\beta}(\Delta)\end{matrix}\right):=
\left(\begin{matrix}\partial_\beta f_\alpha-\partial_\alpha f_\beta+D_i\partial^i_\alpha f_\beta-D_iD_j\partial^{ij}_\alpha f_\beta\\
\partial^i_\beta f_\alpha+\partial^i_\alpha f_\beta-2D_j\partial^{ij}_\alpha f_\beta\\
\partial^{ij}_\beta f_\alpha-\partial^{ij}_\alpha f_\beta\end{matrix}\right).\label{45}
\end{align}
Notice that the expressions \eqref{45} are also defined when $(f_\alpha)$ is not variational. The Helmholtz expressions, as well as the Helmholtz conditions, can be extended to arbitrary order, for example, see \cite{bib36,bib12}.\\ 
\hspace*{0.5cm} Since the the Helmholtz conditions (and expressions) will be very important later on, we give a brief overview how they can be derived and what the main idea is (the idea is pretty much the same as one derives integrability conditions for exact differential forms). According to \eqref{27}, we get that when $I$ is a zero-functional, then $\mathcal{L}_{\text{j}^kV}I_{\sigma}=K_\sigma(\text{j}^kV)$ is again a zero-functional for every $\pi$-vertical vector field $V\in\mathfrak{X}(E)$ with $\operatorname{supp}V\subset\pi^{-1}\Omega$, and in \eqref{andin2} we defined how to apply the formal Lie derivative to such functionals. By direct computation it then follows that the equation
\begin{align}
\mathcal{L}_{\text{j}^kW}\mathcal{L}_{\text{j}^kV}I_{\sigma}-
\mathcal{L}_{\text{j}^kV}\mathcal{L}_{\text{j}^kW}I_{\sigma}-
\mathcal{L}_{\text{j}^k[W,V]}I_{\sigma}=0,\label{33}
\end{align}
must be satisfied for every $\pi$-vertical vector fields $V,W\in\mathfrak{X}(E)$ with $\operatorname{supp}V,W\subset\pi^{-1}\Omega$ and $[W,V]$ is the Lie bracket of $V,W$. Every zero-functional $I$ satisfies \eqref{33}, but not every one-functional $K$ satisfies the analogous equation
\begin{align}
\mathcal{L}_{\text{j}^kW}(K_{\sigma}(\text{j}^kV))-
\mathcal{L}_{\text{j}^kV}(K_{\sigma}(\text{j}^kW))-
K_{\sigma}(\text{j}^k[W,V])=0,\label{34}
\end{align}
for all $\pi$-vertical $V,W\in\mathfrak{X}(E)$ with $\operatorname{supp}V,W\subset\pi^{-1}\Omega$ (since this is only the case when $K_\sigma(\text{j}^kV)=\mathcal{L}_{\text{j}^kV}I_\sigma$, at least locally). Therefore, the condition \eqref{34} delivers a necessary condition for a source form to be variational (it can be shown that it also delivers a locally sufficient condition). With the help of \eqref{34}, one then deduces the Helmholtz conditions with the following procedure. First, \eqref{34} can be written as
\begin{align}
\int_\Omega \text{j}^k\sigma^\ast[(V^\alpha\mathcal{L}_{\text{j}^kW}f_\alpha-W^\alpha\mathcal{L}_{\text{j}^kV}f_\alpha)dx]=0\label{35}
\end{align}
and the integrand in \eqref{35} can be written as
\begin{align}
V^\alpha\mathcal{L}_{\text{j}^kW}f_\alpha-W^\alpha\mathcal{L}_{\text{j}^kV}f_\alpha=&V^\alpha[W^\beta\partial_\beta+(D_iW^\beta)\partial^i_\beta+(D_iD_jW^\beta)\partial^{ij}_\beta]f_\alpha-\nonumber\\
-&W^\alpha[V^\beta\partial_\beta+(D_iV^\beta)\partial^i_\beta+(D_iD_jV^\beta)\partial^{ij}_\beta]f_\alpha.\label{awayof}
\end{align}
Second, the main technique is then to do partial integration with the expression \eqref{awayof} in a systematic way, namely we shift all total derivatives of $D_iV^\beta$ and $D_iD_jV^\beta$ to other terms plus a divergence term (in a similar way as one does partial integration in the first variational formula to derive the Euler-Lagrange equations, where one shifts all total derivatives away of the test functions). It turns out that we get
\begin{align}
&V^\alpha\mathcal{L}_{\text{j}^kW}f_\alpha-W^\alpha\mathcal{L}_{\text{j}^kV}f_\alpha=\nonumber\\
=&V^\alpha W^\beta H_{\alpha\beta}+V^\alpha (D_iW^\beta) H^i_{\alpha\beta}+V^\alpha (D_jD_iW^\beta) H^{ij}_{\alpha\beta}+\nonumber\\
+&D_i[V^\beta D_j(W^\alpha\partial^{ij}_\beta f_\alpha)-W^\alpha V^\beta\partial^i_\beta f_\alpha-W^\alpha (D_jV^\beta)\partial^{ij}_\beta f_\alpha].\label{n56}
\end{align}
Equation \eqref{n56} can be seen as the definition of the Helmholtz expressions. Now we evaluate the integral in \eqref{35}. By Gauss\rq s theorem, the divergence expression 
\begin{align*}
D_i[V^\beta D_j(W^\alpha\partial^{ij}_\beta f_\alpha)-W^\alpha V^\beta\partial^i_\beta f_\alpha-W^\alpha (D_jV^\beta)\partial^{ij}_\beta f_\alpha]
\end{align*}
vanishes when integrated, since $V,W$ have support in $\pi^{-1}\Omega$. The remaining three terms in \eqref{n56} include the Helmholtz expressions, and the coefficients $V^\alpha W^\beta$, $V^\alpha D_iW^\beta$ and $V^\alpha D_iD_jW^\beta$. It can be shown that these three coefficients can be chosen independently and therefore the Helmholtz expressions must vanish under the assumption that \eqref{34} holds. To understand this independence, let us consider 
\begin{align*}
(\text{j}^k\sigma)^\ast(V^\alpha W^\beta)=\varphi^\alpha(x)\psi^\beta(x),\\
(\text{j}^k\sigma)^\ast(V^\alpha D_iW^\beta)=\varphi^\alpha(x)\psi^\beta_i(x),\\
(\text{j}^k\sigma)^\ast(V^\alpha D_iD_jW^\beta)=\varphi^\alpha(x)\psi^\beta_{ij}(x),
\end{align*}
where $\varphi^\alpha(x)$ and $\psi^\beta(x)$ can be seen as test functions and derivatives are written as $\psi^\beta_i=\frac{\partial \psi^\beta}{\partial x^i}$ and $\psi^\beta_{ij}=\frac{\partial^2 \psi^\beta}{\partial x^i\partial x^j}$. Since these test functions can be chosen arbitrarily, we choose $(\varphi^\alpha)=(0,...,0,\varphi^\alpha,0,...,0)$ and $(\psi^\beta)=(0,...,0,\psi^\beta,0,...,0)$, that is, they all vanish except for one index $\alpha$ and $\beta$. Moreover, we chose $\operatorname{supp}\varphi^\alpha\subset$ $\operatorname{supp}\psi^\beta$ and $\psi^\beta(x)\equiv 1$ in the support of $\varphi^\alpha$, that is, derivatives $\psi^\beta_i$ and $\psi^\beta_{ij}$ vanish in the support of $\varphi^\alpha$. Equation \eqref{35} then leads to 
\begin{align*}
\int_\Omega\varphi^\alpha(x)(\text{j}^k\sigma)^\ast(H_{\alpha\beta})dx=0,\quad\beta=1,2,...,m
\end{align*}
for all test functions $\varphi^\alpha$ with the above mentioned properties, and for all sections $\sigma$ which finally leads to $H_{\alpha\beta}=0$ in $J^k\Omega$. A similar discussion then can be done to show that $H^i_{\alpha\beta}=0$ and $H^{ij}_{\alpha\beta}=0$ in $J^k\Omega$.\\ 
\hspace*{0.5cm} When the integrability conditions \eqref{34} are satisfied then, by Poincar\'e lemma, we can construct $I$ (at least locally) as
\begin{align}
I_{\phi_1\sigma}&=\int_{0}^{1}K_{\phi_t\sigma}\left(\frac{d}{dt}\text{j}^k\phi_t\right)dt=\nonumber\\
&=\int_{0}^{1}\left(\int_\Omega f_\alpha(x,s_t(x),Ds_t(x),...,D^k s_t(x))\frac{ds^\alpha_t(x)}{dt} dx\right) dt,\label{functional}
\end{align}
where $\phi_t$ is the flow of a $\pi$-vertical vector field $V\in\mathfrak{X}(E)$ such that $\operatorname{supp}V\subset\pi^{-1}\Omega$, and $(x,s_t(x),Ds_t(x),...,Ds_t(x))$ describes the local coordinates in $J^kE$ of the prolonged section $\phi_t\sigma$ at time $t$.\footnote{Notice that for $\pi$-vertical vector fields we get $\pi_\ast V=0$ and the flow $\phi^0_t$ of $\pi_\ast V$ on $M$ is the identity map. The prolonged transformation of a section $\sigma$ under $\phi_t$ then simplifies to $\text{j}^k[\phi_t\circ\sigma\circ(\phi^0_t)^{-1}](\phi^0_t)=\text{j}^k(\phi_t\circ\sigma)$. And, supp$V\subset\pi^{-1}\Omega$ implies $\frac{d}{dt}\text{j}^k(\phi_t\circ\sigma)\vert_{\pi^{-1}\partial\Omega}=0$ and we can do partial integration without getting any boundary terms.} Interchanging the integrals in \eqref{functional} and choosing coordinates and $\phi_t$ such that it is easy to compute \eqref{functional} leads to the so-called Vainberg-Tonti Lagrangian (it was also discovered earlier by Volterra) \cite[p.5]{bib9}
\begin{align*}
L=\int_0^1 f_\alpha(x,tu,tDu,...,tD^ku)u^\alpha dt.
\end{align*}
Formula \eqref{functional} is very much the same as one constructs the zero-form $\eta$ in \eqref{noticeb}, where we get
\begin{align*}
\eta_{\phi_1(p_0)}=\int_{0}^{1}\omega_{\phi_t(p_0)}\left(\frac{d\phi}{dt}\right)dt,
\end{align*}
and the flow $\phi_t$ connects two points $\phi_0(p_0)=p_0$ and $\phi_1(p_0)=p_1$ in the manifold where $\omega$ is defined.\\
\hspace*{0.5cm} The functional $I$ in \eqref{functional} is defined for any section $\tilde\sigma$ of $\pi$ in the sense that when we choose a fixed section $\sigma$ such that $\tilde\sigma\vert_{\partial\Omega}=\sigma\vert_{\partial\Omega}$, then we can always find a flow $\phi_t$, with the above properties, such that $\phi_0\circ\sigma=\sigma$ and $\tilde\sigma=\phi_1\circ\sigma$. That is, two sections $\tilde\sigma$ and $\sigma$ can only be connected through $\phi_t$ when they have the same values at the boundary $\partial\Omega$. This means that we get different functionals $I$ for different boundary conditions which is for practical reasons no problem, since the boundary conditions are usually fixed, and for any fixed boundary conditions we can always compute the corresponding functional $I$.

\paragraph{Helmholtz dependencies.} As already mentioned, the Helmholtz conditions \eqref{h} guarantee that a second order source form is locally variational, but it turns out that these conditions are not independent and we would actually need less of them to garantee that a source form is locally variational. The main reason for this is that the partial integration technique for two vector fields, to get from \eqref{awayof} to \eqref{n56}, does not necessarily deliver independent Helmholtz expressions (in contrast to the partial integration which is done in the first variational formula, where only one vector field is involved). The Helmholtz expressions \eqref{45} are dependent and we get the following relations
\begin{align}
&H_{\alpha\beta}+H_{\beta\alpha}-D_iH^i_{\alpha\beta}+D_iD_jH^{ij}_{\alpha\beta}=0,\label{58a}\\
&\quad\quad\;\, H^i_{\alpha\beta}-H^i_{\beta\alpha}-2D_jH^{ij}_{\alpha\beta}=0,\label{58b}\\
&\quad\quad\quad\quad\quad\;\, H^{ij}_{\alpha\beta}+H^{ij}_{\beta\alpha}=0\label{58c}
\end{align}
which are always satisfied, whether $(f_\alpha)$ is variational or not (in contrast to the Helmholtz conditions which are only satisfied when $(f_\alpha)$ is variational). We call these relations the \textbf{Helmholtz dependencies}. They can be proven by a straight forward calculation, when substituting the definition \eqref{45} into \eqref{58a}, \eqref{58b} and \eqref{58c}. Using these relations is crucial in our main proof later, since they reduce drastically the number of unknowns. These relations can be found in \cite[p.377]{bib14} and \cite[p.8]{bib1}, and they are already used in \cite{bib17}, but not explicitly mentioned.\\
\hspace*{5mm} Notice that there are also integrability conditions for functions $G_{\alpha\beta},G^i_{\alpha\beta},G^{ij}_{\alpha\beta}$ under which they can be written as \eqref{45} for some functions $(f_\alpha)$. These conditions have to be distinguished from the Helmholtz dependencies which are not these integrability conditions. See \cite[p.86]{bib18}, where these integrability conditions can be found for $n=1$ and third order (they are pretty complicated in general).\\
\hspace*{0.5cm} Now we briefly want to investigate how many independent Helmholtz expressions we actually have. Counting the indices in the Helmholtz expressions (for second order source forms), we get that there are $m^2+nm^2+\frac{n(n+1)}{2}m^2$ Helmholtz expressions and we used that $H^{ij}_{\alpha\beta}$ is symmetric in $i,j$. Equation \eqref{58c} delivers $N_3:=\frac{n(n+1)}{2}\frac{m(m+1)}{2}$ possible independent equations, since the equation is symmetric in $\alpha,\beta$ (and symmetric in $i,j$). We can write \eqref{58b} as
\begin{align}
(H^i_{\alpha\beta}-D_jH^{ij}_{\alpha\beta})-(H^i_{\beta\alpha}-D_jH^{ij}_{\beta\alpha})=0,\label{accord}
\end{align}
when applying \eqref{58c}. This shows that \eqref{58b} is skewsymmetric in $\alpha,\beta$ and therefore it delivers $N_2:=n\frac{m(m-1)}{2}$ possible independent equations. According to \eqref{accord}, equation \eqref{58a} can be written as
\begin{align*}
H_{\alpha\beta}+H_{\beta\alpha}
-\frac{1}{2}D_i(H^i_{\alpha\beta}-D_jH^{ij}_{\alpha\beta})
-\frac{1}{2}D_i(H^i_{\beta\alpha}-D_jH^{ij}_{\beta\alpha})
=0
\end{align*}
which shows that \eqref{58a} is symmetric in $\alpha,\beta$ and we therefore get $N_1:=\frac{m(m+1)}{2}$ possible independent equations. When we assume that no further relations between the equations in \eqref{58a}, \eqref{58b} and \eqref{58c}, then there are at most
\begin{align}
&m^2+nm^2+\frac{n(n+1)}{2}m^2-N_1-N_2-N_3=\nonumber\\
=&\frac{m(m-1)}{2}+n\frac{m(m+1)}{2}+\frac{n(n+1)}{2}\frac{m(m-1)}{2}\label{number}
\end{align}
independent Helmholtz expressions. It seems that the explicit construction of independent Helmholtz expressions is not completed and usually always the expressions \eqref{45} are used in the literature. We investigated this problem in parts in \cite{bib52}, where $n=m=1$, but for higher order source forms.

\paragraph{Symmetries.} We say that an one-functional \textbf{$K$ satisfies a symmetry}, if there exists a $\pi$-projectable vector field $V\in\mathfrak{X}(E)$ such that \eqref{andin} vanishes for all sections $\sigma$ and for all $W\in\mathfrak{X}(J^kE)$. This is equivalent to $\mathcal{L}_{\text{j}^kV}\Delta=0$ (which is a version of the Noether-Bessel-Hagen equation, see \cite[p.177]{bib36}). Consequently, we define a \textbf{symmetry of $\Delta$} as a $\pi$-projectable vector field $V\in\mathfrak{X}(E)$ such that $\mathcal{L}_{\text{j}^kV}\Delta=0\quad\text{for all points in $J^kE$}$. Notice that this definition is in general not equivalent to $\mathcal{L}_{\text{j}^kV}f_\alpha=0$, $\alpha=1,2,...,m$. Roughly speaking, in Takens\rq{} problem we require that the weak formulation of a differential equation satisfies certain symmetries and not the differential equation itself. And, the symmetry does not only hold for solutions of the differential equation, it holds for all points in $J^kE$.

\paragraph{Continuity equations and the total divergence operator.} In the example in Section 2 we used a very vague concept of conservation law and in this paragraph we will give a precise definition. First of all, it is not obvious how to use the intuitive concept of conservation law in physics to get a meaningful mathematical definition which covers all these physical 'conservation laws'. Moreover, ODEs and PDEs can behave quite differently in this regard and, in our opinion, we should distinguish these cases carefully and not give certain equations the same name conservation laws. In the case of ODEs, a conservation law is usually considered as a first integral of the differential equation, that is, a function which is constant along solutions of the differential equation, see the example in Section 2. For PDEs, a conservation law is usually considered as a divergence expression, or in physics known as a continuity equation, and this does in general not deliver a function which is constant along solutions. A well-known example of a continuity equation in physics is 
\begin{align}
\frac{\partial\rho}{\partial t}+\nabla \textbf{j}=0,\label{e1}
\end{align}
where $\rho$ is the charge-, and $\textbf{j}$ the current density in Maxwell\rq s equations. In relativistic notation, this can also be written as $\partial_\mu J^\mu=0$, where $(J^\mu)=(\rho,\,\textbf{j})$ is the so-called \textit{four-current} in physics.\footnote{In our notation $\partial_\mu$ would be the total derivative operator $D_\mu$.} We will just call it \textit{current density}, see below. Notice that in relativistic notation the time and spatial variables are considered equally and the whole equation \eqref{e1} is just a divergence equation. Another example of a continuity equation is given by \eqref{e1}, where $\rho=\vert\Psi\vert^2$ and $\textbf{j}=\frac{i\hbar}{2m}(\Psi\nabla\Psi^\ast-\Psi^\ast\nabla\Psi)$ in Schrödinger\rq s equation, but there the relativistic formulation is not so easy. These examples, especially the relativistic case, now motivate what we will define as a current density and a continuity equation. For further information we refer to \cite{bib12} and \cite{bib0}.\\
\hspace*{0.5cm} We define a \textbf{current density $J$} as an $(n-1)$-form on $J^kE$, written in local coordinates as
\begin{align}
J=\sum_{i=1}^{n}(-1)^{i+1}J^idx^1\wedge ...\wedge \widehat{dx^i}\wedge ...\wedge dx^n,\label{44}
\end{align}
where the hat $\widehat{dx^i}$ denotes omission, and $J^i$ are functions on $J^kE$, see \cite[p.199]{bib12}. The intrinsic property of \eqref{44} is that $J$ is a so-called \textit{horizontal} $(n-1)$-form on $J^kE$ which means that there are no $du^\alpha, du^\alpha_i,...,du^\alpha_{i_1...i_k}$-forms involved, see \cite[p.33]{bib9} and \cite[p.36]{bib36}. For any horizontal $(n-1)$-form $\omega\in\Omega^{n-1}(J^kE)$, written in local coordinates as $\omega=\omega_{i_1i_2...i_{n-1}}dx^{i_1}\wedge ...\wedge dx^{i_{n-1}}$, we define the \textbf{total divergence operator $E_n$} in local coordinates as
\begin{align*}
E_n\omega:=(D_i\omega_{i_1i_2...i_{n-1}})dx^i\wedge
dx^{i_1}\wedge ...\wedge dx^{i_{n-1}},
\end{align*}
where $D_i$ is the total derivative operator in \eqref{divop}. The fact that $E_n$ is coordinate invariant can be seen when writing $E_n$ as $E_n=hd$, where $h$ is the so-called horizontalization operator and $d$ is the exterior derivative, see \cite[p.289]{bib36}, \cite[p.33]{bib15} and \cite[p.54]{bib36}. For any current density $J$ we get $E_nJ=(D_iJ^i)dx$, that is, where the $dx^j$-forms are ordered in an increasing manner, and $D_iJ^i$ delivers a divergence expression which is coordinate invariantly defined by these forms.
\begin{definition}\label{contequ}
We say that a source form $\Delta$ satisfies a \textbf{local continuity equation}, if for each $p\in E$ there exists a neighborhood $U\subset E$, a current density $J$ on $(\pi^{k,0})^{-1}U\subset J^kE$, and a corresponding $\pi^k$-vertical vector field $Q$ on $(\pi^{k,0})^{-1}U\subset J^kE$, such that the equation $E_nJ=Q\lrcorner\Delta$ is satisfied on $(\pi^{k,0})^{-1}U\subset J^kE$. If there exists a current density $J$ on $J^kE$, and a corresponding $\pi^k$-vertical vector field $Q\in\mathfrak{X}(J^kE)$, such that $E_nJ=Q\lrcorner\Delta$ is satisfied on $J^kE$, then we say that $\Delta$ satisfies a \textbf{global continuity equation}.
\end{definition}
The vector field $Q=Q^\alpha\partial_\alpha+...+Q^\alpha_{i_1...i_k}\partial^{i_1...i_k}_\alpha$ is called \textbf{characteristic vector field for the continuity equation} and $(Q^\alpha)=(Q^\alpha(x,u,Du,...,D^k u))$ is called the \textbf{characteristic of the continuity equation}, see \cite[p.199]{bib12} and \cite[p.270]{bib0}. Notice that the equation $E_nJ=Q\lrcorner\Delta$, or in local coordinates $D_iJ^i=Q^\alpha f_\alpha$, can be understood as follows. If we are on solutions of the (differential) equation $f_\alpha=0$, $\alpha=1,2,...,m$, then we get $D_iJ^i=0$, and this is interpreted as a continuity equation which leads to some kind of conserved quantity, like charge conservation in Maxwell\rq s-, or Schrödinger\rq s equations.\\
\hspace*{0.5cm} Since $\Delta$ corresponds to a one-functional $K$, we can also define what it means that $K$ satisfies a continuity equation. Let us assume that the equation $E_nJ=Q\lrcorner\Delta$ is satisfied for some $J$ and $Q$. Then, by Gauss theorem, we get
\begin{align*}
K_\sigma(Q)&=\int_\Omega\text{j}^k\sigma^\ast (f_\alpha Q^\alpha dx)
=\int_\Omega\text{j}^k\sigma^\ast (D_iJ^i dx)
\overset{\text{Gauss}}{=}\int_{\partial\Omega}\text{j}^k\sigma^\ast (J^in_idS)=\\
&=K_\sigma(Q)\vert_{\partial\Omega},
\end{align*}
where $(n_i)$ is the unit normal vector of $\partial\Omega$ and $dS$ the surface volume form on $\partial\Omega$. Therefore, we say that a one-functional $K$ satisfies a \textbf{continuity equation}, if there exists a $\pi^k$-vertical vector field $Q\in\mathfrak{X}(J^kE)$ such that $K_\sigma(Q)$ only depends on the restriction of $\partial\Omega$ for all sections $\sigma$ of $\pi$.\footnote{This also implies that $I_\sigma=K_\sigma(Q)$ is a trivial zero-functional with trivial Euler-Lagrange equation and $L=D_iJ^i$ is a trivial Lagrangian. This will also be discussed below.}

\paragraph{The Euler-Lagrange and Helmholtz operator.} Let us define the \textbf{Euler-Lagrange operator $\mathcal{E}_\alpha$} in local coordinates as
\begin{align*}
\mathcal{E}_\alpha:=\partial_\alpha-D_{i_1}\partial^{i_1}_\alpha+...
+(-1)^kD_{i_1}D_{i_2}...D_{i_k}\partial^{i_1i_2...i_k}_\alpha,\quad \alpha=1,2,...,m,
\end{align*}
which acts on Lagrangians $L=L(x,u,...,D^ku)$ as $\mathcal{E}_\alpha L$. To define the operator in a coordinate invariant way, we define the operator $E_{n+1}$, acting on Lagrange forms $\lambda=Ldx$ as 
\begin{align}
E_{n+1}(\lambda):=(\mathcal{E}_\alpha L)du^\alpha\wedge dx.\label{euler}
\end{align}
Let us also define the \textbf{Helmholtz operator $\mathcal{H}^\gamma_{\alpha\beta}$} (for second order source forms) in local coordinates as
\begin{align*}
\mathcal{H}^\gamma_{\alpha\beta}:=
\left(\begin{matrix}
\mathcal{H}^\gamma_{\alpha\beta}\\\mathcal{H}^{\gamma,i}_{\alpha\beta}\\\mathcal{H}^{\gamma,ij}_{\alpha\beta}\\
\end{matrix}\right)=
\left(\begin{matrix}
\partial_\beta\delta^\gamma_\alpha-\partial_\alpha\delta^\gamma_\beta+
D_{i_1}\partial^{i_1}_\alpha\delta^\gamma_\beta-
D_{i_1}D_{i_2}\partial^{i_1i_2}_\alpha\delta^\gamma_\beta\\
\partial^i_\beta\delta^\gamma_\alpha+\partial^i_\alpha\delta^\gamma_\beta-2
D_{i_1}\partial^{ii_1}_\alpha\delta^\gamma_\beta\\
\partial^{ij}_\beta\delta^\gamma_\alpha-\partial^{ij}_\alpha\delta^\gamma_\beta
\end{matrix}\right),
\end{align*}
where $\alpha,\beta,\gamma=1,2,...,m$ which acts on functions $(f_\gamma)$ as $\mathcal{H}^\gamma_{\alpha\beta} f_\gamma$, such that we get the Helmholtz expressions in \eqref{45} ($\delta^\gamma_\alpha$ denotes the Kronecker-delta). Notice that the operator $\mathcal{H}^\gamma_{\alpha\beta}$ can be defined for higher order $(f_\gamma)$ when following the ideas where we derived the Helmholtz conditions. Again, to define the operator in a coordinate invariant way, we define the operator $E_{n+2}$ acting on (second order) source forms $\Delta=f_\alpha du^\alpha\wedge dx$ as 
\begin{align}
E_{n+2}(\Delta):=&\frac{1}{2}[(\mathcal{H}^\gamma_{\alpha \beta}f_\gamma)du^\beta+
(\mathcal{H}^{\gamma,i}_{\alpha \beta}f_\gamma)du^\beta_i
+(\mathcal{H}^{\gamma,ij}_{\alpha \beta}f_\gamma)du^\beta_{ij}]\wedge du^\alpha\wedge dx.\label{helm}
\end{align}
The fact that $E_{n+1}$ and $E_{n+2}$ are coordinate invariant can be shown when writing $E_{n+1}=\mathcal{I}d$, acting on certain kinds of differential forms, where $\mathcal{I}$ is the so-called interior Euler-Lagrange operator and $d$ is the exterior derivative. The same holds for $E_{n+2}$. See \cite{bib32} and \cite[p.7]{bib30} for the finite jet bundle $J^kE$, and \cite{bib7} for the infinite jet bundle $J^\infty E$. The operator $\mathcal{I}$ is basically needed to do partial integration, as we explained, for example, in \eqref{n56} when we derived the Helmholtz expressions.

\paragraph{Locally exact sequences.} It is clear that $\mathcal{H}^\gamma_{\alpha\beta}\mathcal{E}_\gamma L=0$ for every Lagrangian $L$, since variational expressions $f_\gamma=\mathcal{E}_\gamma L$ satisfy the Helmholtz conditions, see \eqref{h}. As an operator identity we therefore get $\mathcal{H}^\gamma_{\alpha\beta}\mathcal{E}_\gamma=0$, or $E_{n+2}\circ E_{n+1}=0$. It can also be shown that $\mathcal{E}_\alpha D_i=0$, or $E_{n+1}\circ E_n=0$. This is because any Lagrangian of the form $L=D_iJ^i$, where $(J^i)$ are some functions on $J^kE$, leads to a zero-functional of the form
\begin{align*}
I_\sigma=\int_\Omega(\text{j}^k\sigma)^\ast [(D_iJ^i)dx]\overset{\text{Gauss}}{=}\oint_{\partial\Omega}(\text{j}^k\sigma)^\ast (J^i n_idS)=I_\sigma \vert_{\partial\Omega},
\end{align*}
where, by Gauss\rq{}s theorem, $I_\sigma$ does only depend on the boundary of $\Omega$ and therefore perturbations of $I_\sigma\vert_{\partial\Omega}$, which are described by the Euler-Lagrange operator $\mathcal{E}_\alpha$ and $\pi$-vertical vector fields which have support in $\pi^{-1}\Omega$, lead to the trivial Euler-Lagrange equation $\mathcal{E}_\alpha (D_iJ^i)=0$ for any functions $(J^i)$ (notice $n_i$ are the components of the unit normal vector of the surface $\partial\Omega$ and $dS$ is the volume form on $\partial\Omega$).\footnote{To prove $\mathcal{E}_\alpha D_i=0$ (and also to prove $\mathcal{H}^\gamma_{\alpha\beta}\mathcal{E}_\gamma=0$) it is also possible to commute $D_i$ with the partial derivatives $\partial^{i_1...i_l}_\alpha$ in the operator $\mathcal{E}_\alpha$ (or in $\mathcal{H}^\gamma_{\alpha\beta}$) and use the commutation relation $[\partial^{i_1...i_l}_\alpha,D_i]=\delta^\beta_\alpha\delta^{i_1...i_l}_{ij_1...j_{l-1}}\partial^{j_1...j_{l-1}}_\beta$. However, we usually prefer to prove these identities with the methods explained above, since commuting these operators can get very complicated, at least to show that $\mathcal{H}^\gamma_{\alpha\beta}\mathcal{E}_\alpha=0$.} Moreover, the operators $E_n$, $E_{n+1}$ and $E_{n+2}$ lead to a locally exact \textbf{sequence of differential forms}\footnote{To construct $J$, see the total homotopy operator in \cite[p.364]{bib0} in equation (5.112). To construct $\lambda$, we can use the so-called Vainberg-Tonti Lagrangian, see \cite[p.56]{bib9}, \cite[p.136]{bib36} and also see formula \eqref{functional}.} of the form
\begin{align}
 ...\longrightarrow
\{ J\}\overset{E_n}{\longrightarrow}\{\lambda\}\overset{E_{n+1}}{\longrightarrow}\{\Delta\}
\overset{E_{n+2}}{\longrightarrow}...,\label{varseq}
\end{align}
where $\{J\}$ denotes the set of current densities, $\{\lambda\}$ the set of Lagrange forms, and $\{\Delta\}$ the set of source forms. Notice that the operators $E_n$, $E_{n+1}$ and $E_{n+2}$ can be extended to more general differential forms on $J^kE$, see \cite{bib36}. There are actually different ways how to realize this sequence, or sequences similar to them. For example, Krupka uses a finite order sequence with quotient spaces \cite{bib31,bib36}, and Anderson uses a infinite order sequence which is constructed with the help of vertical and horizontal exterior derivatives \cite{bib7}.\\
\hspace*{0.5cm} We can also construct an analog \textbf{sequence of functionals}. Let us define the operator $\tilde\delta_1$ for zero-functionals as
\begin{align}
\tilde\delta_1I_\sigma(\text{j}^kV):=\mathcal{L}_{\text{j}^kV}I_\sigma\label{u1}
\end{align}
for all sections $\sigma$ of $\pi$, and for all $\pi$-vertical vector fields $V\in\mathfrak{X}(E)$ with supp$V\subset\pi^{-1}\Omega$. Moreover, we define $\tilde\delta_2$ for one-functionals as
\begin{align}
\tilde\delta_2K(\text{j}^kV,\text{j}^kW):=
\mathcal{L}_{\text{j}^kV}(K_{\sigma}(\text{j}^kW))-
\mathcal{L}_{\text{j}^kW}(K_{\sigma}(\text{j}^kV))-
K_{\sigma}(\text{j}^k[V,W])\label{u2}
\end{align}
for all sections $\sigma$ of $\pi$, and for all $\pi$-vertical vector fields $V,W\in\mathfrak{X}(E)$ with supp$V,W\subset\pi^{-1}\Omega$. Then we get a locally exact sequence of the form
\begin{align}
 ...\longrightarrow
\{ ?\}\overset{\tilde\delta_0}{\longrightarrow}\{I\}\overset{\tilde\delta_1}{\longrightarrow}\{K\}
\overset{\tilde\delta_2}{\longrightarrow}...,\label{varseq2}
\end{align}
where $\{I\}$ is the set of zero-functionals, $\{K\}$ is the set of one-functionals and so on. The operators $\tilde\delta_k$, where $k\ge 3$, can be defined in a similar way as \eqref{u1} and \eqref{u2} for $\pi$-vertical vector fields with support in $\pi^{-1}\Omega$ and there are corresponding operators $E_{n+k}$ (this is actually the motivation how to define $E_{n+k}$ for $k\ge 1$). However, it turns out that $E_n$ has no direct analog $\tilde\delta_0$ and the sequence \eqref{varseq2} fails here somehow. This is one reason why it can be reasonable to consider the sequence \eqref{varseq} of differential forms instead of the sequence \eqref{varseq2} of functionals. On the other hand, for the definition of a Cartan\rq s formula it is better to consider an \textbf{extended sequence of functionals}, as we will see in the next paragraph, and the sequence of differential forms fails here (Cartan\rq s formula also works for the operators $\tilde\delta_k$, but we need more general operators in Takens\rq{} problem later on). In the extended sequence we will write $\delta_k$ instead of $\tilde\delta_k$, since the operators will be define more generally for $\pi$-projectable vector fields instead of $\pi$-vertical vector fields and we also do not have the restriction that the support must be in $\pi^{-1}\Omega$. Notice that when $V,W\in\mathfrak{X}(E)$ are $\pi$-vertical then $[V,W]$ is also a $\pi$-vertical vector field, and if supp$V,W\subset\pi^{-1}\Omega$ then we also get supp$[V,W]\subset\pi^{-1}\Omega$. This means that the subset of vector fields we consider to define $\tilde\delta_k$ is closed under the Lie bracket and our definition is meaningful.

\paragraph{Cartan\rq s formula and its generalizations.} As it is well-known, by Cartan\rq s formula, the Lie derivative decomposes into two terms $\mathcal{L}_V=d V\lrcorner+V\lrcorner d$. We want to find a similar formula with the formal Lie derivative \eqref{andin} applied to one-functionals. To accomplish this we define a formal exterior derivative $\delta_1$ for zero-functional as
\begin{align}
\delta_1I_\sigma(\text{j}^kV):=\mathcal{L}_{\text{j}^kV}I_\sigma\label{usual1}
\end{align}
for all sections $\sigma$ of $\pi$, and for all $\pi$-projectable vector fields $V\in\mathfrak{X}(E)$. Moreover, we define $\delta_2$ for one-functional as
\begin{align}
\delta_2K_\sigma(\text{j}^kV,\text{j}^kW):=
\mathcal{L}_{\text{j}^kV}(K_{\sigma}(\text{j}^kW))-
\mathcal{L}_{\text{j}^kW}(K_{\sigma}(\text{j}^kV))-
K_{\sigma}(\text{j}^k[V,W])\label{usual2}
\end{align}
for all sections $\sigma$ of $\pi$, and for all $\pi$-projectable vector fields $W,V\in\mathfrak{X}(E)$. It holds $\delta_2\delta_1=0$. Now we get a formal Cartan\rq s formula
\begin{align}
\mathcal{L}_{\text{j}^kV}=\delta_1(\text{j}^kV)\lrcorner+(\text{j}^kV)\lrcorner\delta_2\label{lack3}
\end{align}
which can be applied to one-functionals. To prove this, we can write
\begin{align}
&(\mathcal{L}_{\text{j}^kV}K_\sigma)(\text{j}^kW)
=(\delta_1K_\sigma(\text{j}^kV))(\text{j}^kW)+\delta_2K_\sigma(\text{j}^kV,\text{j}^kW)=\label{cartan2}\\
=&\mathcal{L}_{\text{j}^kW}(K_\sigma(\text{j}^kV))
+\mathcal{L}_{\text{j}^kV}(K_{\sigma}(\text{j}^kW))-
\mathcal{L}_{\text{j}^kW}(K_{\sigma}(\text{j}^kV))-
K_{\sigma}(\text{j}^k[V,W])=\nonumber\\
=&\mathcal{L}_{\text{j}^kV}(K_{\sigma}(\text{j}^kW))-
K_{\sigma}(\text{j}^k[V,W])\nonumber
\end{align}
and one can verify the last of these equations by direct calculation. Notice that there is no such nice formula with the operators $E_{n+1}$ and $E_{n+2}$, that is,
\begin{align}
\mathcal{L}_{\text{j}^kV}\ne E_{n+1}(\text{j}^kV)\lrcorner+(\text{j}^kV)\lrcorner E_{n+2},\label{lack2}
\end{align}
when applied to source forms $\Delta$ (not even when $V$ is $\pi$-vertical and supp$V\subset\pi^{-1}\Omega$). The lack of such formulas with the operators $E_{n+1}$ and $E_{n+2}$ is why we prefer to use $\delta_1$ and $\delta_2$ to explain Takens\rq{} problem. Also notice that there is no such nice formula with the operators $E_n$ and $E_{n+1}$, that is,
\begin{align}
\mathcal{L}_{\text{j}^kV}\ne E_{n}(\text{j}^kV)\lrcorner+(\text{j}^kV)\lrcorner E_{n+1},\label{lack1}
\end{align}
when applied to Lagrange forms $\lambda$. Formula \eqref{lack1} can be fixed in some sense as $h\mathcal{L}_{\text{j}^kV}=h E_{n}(\text{j}^kV)\lrcorner+h(\text{j}^kV)\lrcorner E_{n+1}$, where $h$ is the so-called horizontalization operator, but then one also has to apply it to so-called Lepage forms, see \cite[p.173]{bib36}. As far as we know, there is no formula noticed in the literature which fixes \eqref{lack2}. However, there is a (coordinate dependent) formula which is at least very similar to Cartan\rq s formula and which is given in the next theorem.
\begin{theorem}\label{t45}
Let $\Delta$ be a second order source form, and $V\in\mathfrak{X}(E)$ a $\pi$-projectable vector field. Then we get the decomposition
\begin{align*}
\mathcal{L}_{\text{j}^kV}\Delta
=\{\mathcal{E}_\alpha(V^\beta_{ch} f_\beta)
+[V^\beta_{ch}\mathcal{H}^\gamma_{\alpha\beta}+(D_iV^\beta_{ch})\mathcal{H}^{\gamma,i}_{\alpha\beta}+(D_jD_iV^\beta_{ch})\mathcal{H}^{\gamma,ji}_{\alpha\beta}]f_\gamma\}du^\alpha\wedge dx.
\end{align*}
\end{theorem}
Notice that on the left hand side of the equation in Theorem \ref{t45} the vector field $\text{j}^kV$ occurs, and on the right hand side only $\text{j}^kV_{ch}$ occurs. This also shows that that a Cartan\rq s formula with these differential forms, if it exists, is in any case more complicated than the usual one. Theorem \ref{t45} can be found in \cite[p.202]{bib12} and it also holds for higher order source forms.\footnote{Notice that $H_\Delta(Y)$ in \cite[p.201]{bib12} may be given coordinate invariantly, but there is no corresponding Cartan\rq s formula with operators $E_{n+1}$ and $E_{n+2}$.} It also delivers an implicit formulation of Noether\rq{}s theorem.

\paragraph{The formulation of Takens\rq{} problem.} Now let us explain Takens\rq{} problem \cite{bib17}. The key to understand it is Cartan\rq s formula \eqref{cartan2}, where the vector fields $V$ and $W$ are of the following form
\begin{enumerate}
\item $V\in\mathfrak{X}(E)$ is a $\pi$-projectable symmetrie vector field.
\item The equation holds for all $\pi$-vertical $W\in\mathfrak{X}(E)$ with supp$W\subset\pi^{-1}\Omega$.
\end{enumerate}
Roughly speaking, we assume some kinds of symmetries and continuity equations for a one-functional $K$, and we want to prove that $K$ is variational. By definition, a symmetry of $K$ is an equation $\mathcal{L}_{\text{j}^kV}K_\sigma=0$ which means that the left hand side in \eqref{cartan2} vanishes for all $\text{j}^kW$ (even if we require the equation only for $\text{j}^kW$, where $W$ is $\pi$-vertical and supp$W\subset\pi^{-1}\Omega$). Next, we consider the expression $(\delta_1 K_\sigma(\text{j}^kV))(\text{j}^kW)$ on the right hand side of \eqref{cartan2}. There we can write
\begin{align*}
K_\sigma(\text{j}^kV)
&=\int_\Omega(\text{j}^k\sigma)^\ast [(\text{j}^kV)\lrcorner\Delta]=\\
&=\int_\Omega(\text{j}^k\sigma)^\ast[V^\alpha f_\alpha dx+(-1)^iV^if_\alpha du^\alpha\wedge dx^1\wedge ...\widehat{dx^i}...\wedge dx^n]=\\
&=\int_\Omega(\text{j}^k\sigma)^\ast[V^\alpha f_\alpha dx+(-1)^iV^if_\alpha u^\alpha_idx^i\wedge dx^1\wedge ...\widehat{dx^i}...\wedge dx^n]=\\
&=\int_\Omega(\text{j}^k\sigma)^\ast[(V^\alpha-V^iu^\alpha_i)f_\alpha dx]=\\
&=K_\sigma(\text{j}^kV_{ch})
\end{align*}
and if $\text{j}^kV_{ch}=Q$ is the characteristic vector field for a continuity equation of $K$, then $K_\sigma(\text{j}^kV_{ch})=K_\sigma(\text{j}^kV_{ch})\vert_{\partial\Omega}$, that is, $K_\sigma(\text{j}^kV_{ch})$ only depends on $\partial\Omega$.\footnote{Notice that $\text{j}^kV_{ch}$ is a $\pi^k$-vertical vector field, as required in our Definition of a continuity equation.} Therefore, $K_\sigma(\text{j}^kV)$ is a trivial zero-functional such that the first variation vanishes which means that $(\delta_1K_\sigma(\text{j}^kV))(\text{j}^kW)=0$ in \eqref{cartan2} for all $\pi$-vertical $W\in\mathfrak{X}(E)$ with supp$W\subset\pi^{-1}W$ (here we need the additional conditions for $W$). Thus, by symmetry and continuity equation assumptions we get the equation
\begin{align}
\delta_2K_\sigma(\text{j}^kV,\text{j}^kW)=0\label{deltaK}
\end{align}
for all $\pi$-projectable symmetries $V$ and for all $\pi$-vertical $W\in\mathfrak{X}(E)$ with supp$W\subset\pi^{-1}\Omega$. Now the idea is that for sufficiently many symmetries $V$ we can conclude from \eqref{deltaK} that
\begin{align*}
\tilde\delta_2K=0
\end{align*}
which means that $K$ is locally variational. Of course, $\tilde\delta_2K=0$ and $\delta_2K=0$ are not equivalent, but when $\delta_2K_\sigma(Z,\text{j}^kW)=0$ for all $\pi$-vertical $W\in\mathfrak{X}(E)$ with supp$W\subset\pi^{-1}\Omega$, and for a large set of vector fields $Z$ (for example, for all vector fields $Z\in\mathfrak{X}(J^kE)$), then we may also be able to prove that the (restricted) equations 
\begin{align*}
\text{$\delta_2K_\sigma(\text{j}^kY,\text{j}^kW)=0$ for all $\pi$-vertical $Y,W\in\mathfrak{X}(E)$ with supp$Y,W\subset\pi^{-1}\Omega$}
\end{align*}
are satisfied which is the same as $\tilde\delta_2K=0$. Especially the fact that the symmetries $V$ in equation \eqref{deltaK} can have support outside of $\pi^{-1}\Omega$ is no major issue, since outside the support of $W$ the values of $V$ do not matter and for any fixed $W$ we can choose an equivalent $V$ with support in $\pi^{-1}\Omega$ such that we get the same equation \eqref{deltaK}. Also the fact that $V$ are $\pi$-projectable and not $\pi$-vertical vector fields may not have a major impact on what we want to prove, as long as we have sufficiently many symmetry vector fields and corresponding equations \eqref{deltaK} (in some cases, for example, when we have sufficiently many $\pi$-projectable symmetry vector fields then we may be able to take linear combinations such that we get sufficiently many $\pi$-vertical vector fields).\\
\hspace*{0.5cm} One should be aware that when the symmetry vector field $V$ generates a continuity equation, then the corresponding characteristic $Q^\alpha=V^\alpha_{ch}=V^\alpha-u^\alpha_iV^i$ must have a very special form and the connection between symmetries and continuity equations in Theorem \ref{t47} is a strong assumption.\\
\hspace*{0.5cm} We can the same kind of formulation of Takens\rq{} problem with the source forms. Let us assume that:
\begin{itemize}
\item $\Delta$ satisfies the symmetry $\mathcal{L}_{\text{j}^kV}\Delta=0$ for every point on $J^kE$.
\item For every symmetry $V\in\mathfrak{X}(E)$ we get that $\Delta$ satisfies a corresponding global continuity equation of the form $E_nJ=(\text{j}^kV_{ch})\lrcorner\Delta$, where the characteristic of the continuity equation is $Q^\alpha=V^\alpha_{ch}$. Due to the locally exact sequence, it also works with a local continuity equation such that $E_{n+1}(\text{j}^kV_{ch}\lrcorner\Delta)=\mathcal{E_\alpha}(V^\beta_{ch}f_\beta)du^\alpha\wedge dx=0$. In any case we assume that $\mathcal{E_\alpha}(V^\beta_{ch}f_\beta)=0$.
\item $\Delta$ is of second order (in what we investigate in this paper).
\end{itemize}
Then Theorem \ref{t45} forces the equation
\begin{align}
V^\beta_{ch} H_{\alpha\beta}+(D_iV^\beta_{ch})H^i_{\alpha\beta}+(D_jD_iV^\beta_{ch})H^{ji}_{\alpha\beta}=0.\label{46}
\end{align}
We call \eqref{46} the \textbf{equation of continuities and symmetries (ECS)}. If we can conclude from \eqref{46} that $H_{\alpha\beta},H^i_{\alpha\beta},H^{ij}_{\alpha\beta}=0$, then, due to the locally exact sequence, $\Delta$ must be locally variational.\\ 
\hspace*{0.5cm} In the next section we have different kinds of symmetry vector fields $V_{\mathscr{A}}$, where $\mathscr{A}=1,2,...,R$ label these symmetries. For every symmetry vector field we get an equation \eqref{46}. For sufficiently many symmetries, such that the matrix
\begin{align}
(V^\alpha_{ch,\mathscr{A}},D_lV^\beta_{ch,\mathscr{A}},D_jD_iV^\gamma_{ch,\mathscr{A}})^{\alpha,\beta,\gamma=1,2,...,m}_{l,i,j=1,2,...,n;\mathscr{A}=1,2,...,R}\label{47}
\end{align} 
is invertible on $J^kE$ (as a quadratic matrix), the linear equation \eqref{46} immediately forces that $H_{\alpha\beta},H^i_{\alpha\beta},H^{ij}_{\alpha\beta}=0$ (notice that $l,i,j,\alpha,\beta,\gamma$ label the columns of the matrix \eqref{47} and $\mathscr{A}$ the rows). However, \eqref{47} does in general not have full rank. One reason for this is that the entries of the matrix \eqref{47} are always related by each other through total derivatives $D_l$ and $D_jD_i$ and this means that the columns are in some sense dependent. And, in our main theorem later, the integer $R$ times the free index $\alpha$ in \eqref{46} is much smaller than the number of Helmholtz expressions which also means that the matrix \eqref{47} cannot have full rank. Even when we would try to use non-dependent Helmholtz expressions, then the number of non-dependent Helmholtz expressions would still be much larger than $R$ times the free index $\alpha$, see \eqref{number}. This means a much deeper investigation of \eqref{46} is necessary which makes it a really hard problem.\\
\hspace*{0.5cm} Also showing that $\delta_2K_\sigma(\text{j}^kV,\text{j}^kW)=0$ in \eqref{deltaK} is equivalent to $\tilde\delta_2K_\sigma$ is not trivial, since the set of symmetries $V$ can lead to a very specific set of jet prolongations $\text{j}^kV$ (lets say if the symmetry Lie algebra is finite dimensional), whereas the (infinite dimensional) Lie algebra of arbitrarily $\pi$-vertical vector fields on $E$ could lead to a much larger set of possible jet prolongations on $TJ^kE$. There are also counter examples of third order source forms, where our main Theorem \ref{t47} is no longer true, and this also shows that a much deeper investigation of Takens\rq{} problem is necessary.

\section{The main result}
\label{sec:3}
In this section we formulate and prove our main theorem. To specify our assumptions, let us define the space $\mathcal{V}$ of symmetry vector fields on $E$ as
\begin{align*}
\mathcal{V}:=\{V\in\mathfrak{X}(E): V\text{ is $\pi$-projectable and $\mathcal{L}_{\text{j}^kV}\Delta=0$ on $J^kE$}\}.
\end{align*}
Moreover, let us define the condition
\begin{align}
\text{span}_{\mathbb{R}}\{V_p:V\in\mathcal{V}\}=T_pE\quad\text{for all $p\in E$}.\label{48}
\end{align} 
The main result of this paper is the following.
\begin{theorem}\label{t47}
Let $\pi:E\to M$ be a fiber bundle with base dimension $n$ and fiber dimension $m$, where $n,m\in\mathbb{N}$ are arbitrary. Furthermore, let $\Delta$ be a second order source form on $J^2E$. Assume:
\begin{enumerate}
\item The set $\mathcal{V}$ of symmetries of $\Delta$ satisfies \eqref{48}.
\item For each $V\in\mathcal{V}$ we have a corresponding local continuity equation, that is, $E_{n+1}(\text{j}^kV_{ch}\lrcorner\Delta)=0$, where $V_{ch}$ is the characteristic vector field of $V$.
\end{enumerate}
Then $\Delta$ must be locally variational.
\end{theorem}

\noindent The proof of Theorem \ref{t47} will be mainly done in local coordinates and we have to prepare two propositions and some notation first, before we go to the main part of the proof. Roughly speaking, the set $\mathcal{V}$ of symmetry vector fields can be very large and complicated and the first propositions shows that near a point $p\in E$ it is sufficient to consider a finite set of symmetry vector fields. Moreover, the second proposition shows that for such a finite set of symmetry vector fields it is possible to simplify equation \eqref{46}, that is, we can invert a part of the matrix \eqref{47}.
\begin{proposition}\label{p49}
Let $\mathcal{W}\subset\mathfrak{X}(E)$ be a set of $\pi$-projectable vector fields on $E$ such that $\text{span}_\mathbb{R}\{V_p:V\in\mathcal{W}\}=T_pE$ for each $p\in E$. Then for every $p_0\in E$ there exists a small neighborhood $U_{p_0}\subset E$ of $p_0$ such that we can choose $n+m$ vector fields $\{V_1,V_2,...,V_{n+m}\}\subset\mathcal{W}$ such that $\text{span}_\mathbb{R}\{V_{1,p},V_{2,p},...,V_{n+m,p}\}=T_pE$ for all $p\in U_{p_0}$.
\end{proposition}
The proof is straight forward. In the following, we always label finitely many symmetry vector fields by $\mathscr{A}$ and since we are interested in vector fields which span $T_pE$ at each $p\in E$ we may assume that $\mathscr{A}=1,2,...,n+m$.
\begin{proposition}\label{p50}
Let $V_\mathscr{A}$ be $\pi$-projectable vector fields on $E$ and $\mathscr{A}=1,2,...,n+m$. Furthermore let $B:=(V^i_\mathscr{A},V^\alpha_\mathscr{A})^{i;\alpha}_\mathscr{A}$ be the $(n+m)\times (n+m)$-matrix of coefficients of $V_\mathscr{A}$, where $\mathscr{A}$ labels the rows and $i;\alpha$ label the columns. If $\text{span}_{\mathbb{R}}\{V_{\mathscr{A},p}\;,\;\mathscr{A}=1,2,...,n+m\}=T_pE$ for each $p\in U\subset E$,
then there exists a $(n+m)\times (n+m)$-matrix $C=C(x,u)$ on $U\subset E$, such that the $(n+m)\times (n+m)$-matrix
$B=B(x,u)$ satisfies $C\cdot B=Id$, where Id is the identity matrix. That is, $C$ is the inverse matrix of $B$. Moreover, there exists a row $c^\mathscr{A}=c^\mathscr{A}(x,u)$ of the matrix $C$ such that either $c^\mathscr{A}V^i_\mathscr{A}=\delta^{ij}$, $c^\mathscr{A}V^\alpha_\mathscr{A}=0$, or such that $c^\mathscr{A}V^i_\mathscr{A}=0$, $c^\mathscr{A}V^\alpha_\mathscr{A}=\delta^{\alpha\beta}$, where $\delta^{ij}$ and $\delta^{\alpha\beta}$ are Kronecker deltas.
\end{proposition}
The proof follows directly by definition of $\text{span}_{\mathbb{R}}\{V_{p,\mathscr{A}}\;,\;\mathscr{A}=1,2,...,n+m\}=T_pE$ for each $p\in U\subset E$.

\paragraph{The notation of order.} Varying terms in equations independently from others, and thereby showing that they must vanish, will be one of our main techniques to prove Theorem \ref{t47}. With varying terms independently we mean that the coordinates $x^i,u^\alpha,u^\alpha_{i_1},...,u^\alpha_{i_1i_2...i_k}$ can be varied independently if an equation must hold on $J^kE$ (for example, the equation $\mathcal{L}_{\text{j}^kV}\Delta =0$ on $J^kE$). Usually we start our discussion with the highest order coordinate $u^\alpha_{i_1i_2...i_k}$ and all the lower order coordinates are not relevant in this moment. Therefore, it is reasonable to have a notation for such lower order-, or in other words non-important terms, and we will write $O(k-1)$ which stands for a function of order $k-1$, see \eqref{order}.\\ 
\hspace*{0.5cm} Total derivatives increase the order by one when applied to functions. More precisely, the order is increased affine linear in the highest order coordinates and we get
\begin{align}
D_ig&=\partial_ig+u^\alpha_i\partial_\alpha g+...+u^\alpha_{ii_1i_2...i_k}\partial^{i_1i_2...i_k}_\alpha g=O(k)+u^\alpha_{ii_1i_2...i_k}O(k)\label{51}
\end{align}
for every function $g\in C^\infty(J^kE)$. Since $J^{k+1}E$ is an affine linear bundle over $J^kE$, this notation is invariant under coordinate transformations. We also want to introduce the notation $g= O_1(k)$ if $g$ is affine linear in the $k$-th order coordinates, or in other words a polynomial of degree one in $u^\alpha_{i_1i_2...i_k}$. In general, we write $g= O_P(k)$ if $g$ is a polynomial of degree $P$ in the $k$-th order coordinates. However, notice that later $O_\mathscr{A}(k)$ does not indicate a polynomial of degree $\mathscr{A}$, it rather labels the different kinds of symmetries. Sometimes we will write a few indices on the expression $O(k)$, for example, $O^\beta_{ij}(k)$ or $O^\beta_{\mathscr{A},ij}(k)$, and always when we use the indices $\mathscr{A},\alpha,\beta,i,j$ then we do not describe polynomial degree (this will also be clear from the context). Also notice that the definition of objects $O_P(k)$ is invariant under local coordinate transformations. The $O_P(k)$-notation satisfies the following properties
\begin{align*}
O_{P_1}(k)O_{P_2}(k)&=O_{P_1+P_2}(k),\quad\text{for all $k\ge 1$},\\
O_{P_1}(k)O_{P_2}(l)&=O_{P_1}(k),\quad\text{for all $k>l\ge 0$},\\
D_iO_P(k)&=O_1(k+1).
\end{align*}
\hspace*{0.5cm} Let us now explain where we use this notation. The leading orders and the polynomial degree of local coordinates in each of the rows in the matrix \eqref{47} will be crucial in the proof of Theorem \ref{t47} and we can write
\begin{align*}
V^\beta_{ch}&=V^\beta-u^\beta_{i}V^{i},\\
D_kV^\beta_{ch}&=D_k(V^\beta-u^\beta_{i}V^{i})=O^\beta_k(1)-u^\beta_{ki}V^{i},\\
D_lD_kV^\beta_{ch}&=D_l(O^\beta_k(1)-u^\beta_{ki}V^{i})=
O^\beta_{lk}(2)-u^\beta_{lki}V^{i},
\end{align*}
where $V^{i}=V^{i}(x)$ and $V^\beta=V^\beta(x,u)$.\\
\quad\\
\paragraph{Proof of Theorem \ref{t47}.} Let us consider an arbitrary point $p_0\in E$. According to Proposition \ref{p49}, we can find a small neighborhood $U_{p_0}\subset E$ of $p_0$, such that $n+m$ symmetry vector fields $\{V_\mathscr{A},\mathscr{A}=1,2,...,n+m\}\subset\mathcal{V}$ span $T_pE$ for each $p\in U_{p_0}$, and according to symmetry and continuity equation assumptions, we get equation \eqref{46}, that is,
\begin{align}
(V^\beta_\mathscr{A}-u^\beta_{i}V^{i}_\mathscr{A})H_{\alpha\beta}+
(O^\beta_{i\mathscr{A}}(1)-u^\beta_{ki}V^{i}_\mathscr{A})H^k_{\alpha\beta}+
(O^\beta_{ji\mathscr{A}}(2)-u^\beta_{lki}V^{i}_\mathscr{A})H^{lk}_{\alpha\beta}=0\label{52}
\end{align}
for every $\alpha=1,2,...,m$ and every $\mathscr{A}=1,2,...,n+m$. We will show that equation \eqref{52} forces the Helmholtz expressions to be zero in $(\pi^{k,0})^{-1}U_{p_0}\subset J^kE$ and since we can show this for every open subset $U_{p_0}\subset E$ the Helmholtz conditions must be satisfied everywhere on $(\pi^{k,0})^{-1}E=J^kE$. Then, due to the locally exactness of the variational sequence \eqref{varseq}, we know that $\Delta$ must be locally variational.\\ 
\hspace*{0.5cm} Now we will discuss the ECS \eqref{52} in more detail, and we divide the proof into seven steps. The main results in every step will be written in a box. The results in these boxes will then be needed in the next steps. Things which are not written in boxes are basically the proofs of what is written in the boxes.\\
\quad\\
\underline{Step 1 (transform the ECS):}\footnote{This step was developed during my stay at Utah State University when working together with Ian M. Anderson.} According to Proposition \ref{p50}, we can take linear combinations 
\begin{align}
c^\mathscr{A}(V^\beta_\mathscr{A}-u^\beta_{i}V^{i}_\mathscr{A})H_{\alpha\beta}+
&c^\mathscr{A}(O^\beta_{\mathscr{A},k}(1)-u^\beta_{ki}V^{i}_\mathscr{A})H^k_{\alpha\beta}+\nonumber\\
&\quad\quad\quad\quad\quad\quad\quad\quad+c^\mathscr{A}(O^\beta_{\mathscr{A},kl}(2)-u^\beta_{kli}V^{i}_\mathscr{A})H^{kl}_{\alpha\beta}=0\nonumber
\end{align}
such that we get 
\begin{align*}
\text{i)}\quad 0&=-u^\beta_j H_{\alpha\beta}+(O^\beta_{k}(1)-u^\beta_{jk})H^k_{\alpha\beta}+(O^\beta_{kl}(2)-u^\beta_{jkl})H^{kl}_{\alpha\beta},\\
\text{ii)}\quad 0&=H_{\alpha\gamma}+O^\beta_{\gamma,k}(1)H^k_{\alpha\beta}+O^\beta_{\gamma,kl}(2)H^{kl}_{\alpha\beta}.
\end{align*}
To derive equation i), we use $c^\mathscr{A}V^i_\mathscr{A}=\delta^{ij}$ and $c^\mathscr{A}V^\beta_\mathscr{A}=0$, and to derive equation ii), we use $c^\mathscr{A}V^i_\mathscr{A}=0$ and $c^\mathscr{A}V^\beta_\mathscr{A}=\delta^{\beta\gamma}$. Then we take another linear combination of i) and ii) , namely
\begin{align*}
-&u^\beta_j H_{\alpha\beta}+(O^\beta_{k}(1)-u^\beta_{jk})H^k_{\alpha\beta}+(O^\beta_{kl}(2)-u^\beta_{jkl})H^{kl}_{\alpha\beta}+\\
+&u^\gamma_j(H_{\alpha\gamma}+O^\beta_{\gamma,k}(1)H^k_{\alpha\beta}+O^\beta_{\gamma,kl}(2)H^{kl}_{\alpha\beta})=0,
\end{align*}
to derive the transformed ECS
\begin{empheq}[box=\fbox]{align}
\text{I)}\quad 0&=(O^\beta_{jk}(1)-u^\beta_{jk})H^k_{\alpha\beta}+(O^\beta_{jkl}(2)-u^\beta_{jkl})H^{kl}_{\alpha\beta},\quad\quad\text{for all}\,j,\alpha,\label{53}\\
\text{II)}\quad 0&=H_{\alpha\gamma}+O^\beta_{\gamma,k}(1)H^k_{\alpha\beta}+O^\beta_{\gamma,kl}(2)H^{kl}_{\alpha\beta},\quad\quad\quad\quad\quad\,\,\text{for all}\,\alpha,\gamma,\label{54}
\end{empheq}
where we eliminated $H_{\alpha\beta}$ in I), and furthermore II) is the same as ii). From now on we simply call \eqref{53} equation I) and \eqref{54} equation II) and we should remember these equations since they are used several times. Equations I) and II) are together $m^2+nm$ equations, but $H_{\alpha\beta},H^i_{\alpha\beta},H^{ij}_{\alpha\beta}$ are together $m^2+nm^2+\frac{n(n+1)}{2}\frac{m(m-1)}{2}$ unknowns. Therefore, for large $n,m$ we have a highly under-determined system and only for $m=1$ we can immediately determine the solution. To determine to solution of I) and II) for $m>1$ we have to investigate the deeper structure of $H_{\alpha\beta},H^i_{\alpha\beta},H^{ij}_{\alpha\beta}$ and their relations. Since $(f_\alpha)$ are of second order, the polynomial structure of $H_{\alpha\beta},H^i_{\alpha\beta},H^{ij}_{\alpha\beta}$ in fourth and third order coordinates is given as
\begin{align}
H_{\alpha\beta}&=\partial_\beta f_{\alpha}-\partial_\alpha f_{\beta}+D_k\partial^k_\alpha  f_{\beta}-D_kD_l\partial^{kl}_\alpha f_{\beta}=\nonumber\\
&=O_1(3)-u^\delta_{krs}u^\gamma_{lij}\partial^{rs}_\delta\partial^{ij}_\gamma\partial^{kl}_\alpha f_{\beta}-u^\gamma_{ijkl}\partial^{ij}_\gamma\partial^{kl}_\alpha f_\beta,\label{stern}\\
H^i_{\alpha\beta}&=\partial^i_\beta f_{\alpha}+\partial^i_\alpha f_{\beta}-2D_k\partial^{ik}_\alpha f_{\beta}=O_1(3),\nonumber\\
H^{ij}_{\alpha\beta}&=\partial^{ij}_\beta f_{\alpha}-\partial^{ij}_\alpha f_{\beta}=O(2).\nonumber
\end{align}
In Step 2, we now investigate terms which depend on the fourth order coordinates of degree one and terms which depend on third order coordinates of degree two. These terms only occur in equation II) and they must vanish separately. They are generated by the double total derivatives $D_kD_l$ in $H_{\alpha\beta}$, as we can see in \eqref{stern}. Later on we will also discuss terms which are generated by single total derivatives $D_k$ and after that the remaining terms of second order coordinates. In other words, we discuss successively the leading orders and the leading polynomial degree in equation I) and II).\\
\quad\\
\underline{Step 2 (fourth order of degree one and third order of degree two in II)):}\\
As we already mentioned in Step 1, the fourth order terms of degree one
\begin{align}
u^\gamma_{ijkl}\partial^{ij}_\gamma\partial^{kl}_{\alpha}f_{\beta}=0\quad\Leftrightarrow\quad
\partial^{(ij}_\gamma\partial^{kl)}_\alpha f_\beta=0\label{will}
\end{align}
and the third order terms of degree two
\begin{align}
u^\delta_{krs}u^\gamma_{lij}\partial^{rs}_\delta\partial^{ij}_\gamma\partial^{kl}_\alpha f_{\beta}=0\quad\Leftrightarrow\quad
\partial^{(rs}_\delta\partial^{k)(l}_\alpha\partial^{ij)}_\gamma f_\beta=0\label{will2}
\end{align}
must vanish separately in equation II). The brackets $(...)$ mean symmetrization in the indices $ijkl$, $rsk$ and $lij$. 
It is possible to solve the system of differential equations \eqref{will} and \eqref{will2} and we can determine the most general solution. Surprisingly, the set of solutions can be described by a finite dimensional vector space which is usually not the case for systems of PDEs (so-called overdetermined system). The solution $f_\beta$ satisfies the following properties:
\begin{empheq}[box=\fbox]{align}
&\text{$f_\beta$ must be a polynomial of degree $\le n$ in the second order coordi-}\nonumber\\
&\text{nates and the coefficients of the polynomial are first order func-}\nonumber\\
&\text{tions. That is,}\; 
f_\beta=A_\beta+A^{ij}_{\beta\vert\gamma}u^\gamma_{ij}+...+A^{i_1j_1...i_nj_n}_{\beta\vert\gamma_1...\gamma_n}u^{\gamma_1}_{i_1j_1}...u^{\gamma_{n}}_{i_nj_n},\;\text{where}\nonumber\\
&A^{i_1j_1...i_lj_l}_{\beta\vert\gamma_1...\gamma_l}\;\text{are first order functions.}
\label{60} 
\end{empheq}
There are further restrictions on the first order coefficients of $f_\beta$ and we could also describe this structure in detail, but the additional structure of these coefficients is not needed to complete the proof and we will not bother the reader at this point (the underlying structure is that $f_\beta$ is a sum of so-called Hyperjacobians of second order).\footnote{For $n=1$, equation \eqref{will} reduces to $\partial^{11}_\gamma\partial^{11}_\alpha f_\beta=0$ and an easy integration leads to $f_\beta=A_\beta+B_{\beta\gamma}u^\gamma_{11}$, where $A_\beta,B_{\beta\gamma}$ are first order functions. Therefore, $f_\beta$ is obviously a polynomial of degree one in second order coordinates. Equation \eqref{will2} is automatically satisfied, but this is no longer true for $n>1$ and integrating these differential equations gets much more complicated.} The proof is due to Anderson and Duchamp \cite[pp.786]{bib23}.\\
\quad\\
\underline{Step 3 (third order of degree one in I)):} With the help of equation I) and \eqref{60} we show that
\begin{empheq}[box=\fbox]{align}
&H^i_{\alpha\beta}=O(2),\nonumber\\
&H^{ij}_{\alpha\beta}=0.\label{af}
\end{empheq}
It needs a lot of work to deduce \eqref{af} and therefore we will prove it separately in Section 5.\footnote{This step can maybe also be proven with a modification of the so called $d$-fold operator, used in the proof of Theorem 1.1 in \cite[p.379]{bib14} or in the proof of Theorem 1 in \cite[p.12]{bib16}. The same is maybe true for Step 6.}\\
\quad\\
\underline{Step 4 (third order of degree one in II)):} This step is quite simple. We use \eqref{af} and we plug these expressions into equation II) which then leads to
\begin{align*}
\text{II)}\quad H_{\alpha\gamma}+\underbrace{O^\beta_{\gamma,k}(1)H^k_{\alpha\beta}}_{=O(2),\,\text{see \eqref{af}}}+\underbrace{O^\beta_{\gamma,kl}(2)H^{kl}_{\alpha\beta}}_{=0,\,\text{see \eqref{af}}}=0,
\end{align*}
and therefore we get
\begin{empheq}[box=\fbox]{align}
H_{\alpha\gamma}=O(2).\label{57}
\end{empheq}
\quad\\
\underline{Step 5 (we consider the Helmholtz dependencies \eqref{58a}):} This step is again very simple. We use \eqref{af} and \eqref{57} and we plug them into the Helmholtz dependencies \eqref{58a} which delivers
\begin{align*}
\underbrace{H_{\alpha\beta}+H_{\beta\alpha}}_{=O(2),\,\text{see \eqref{57}}}=D_kH^k_{\alpha\beta}-D_kD_l\underbrace{H^{kl}_{\alpha\beta}}_{=0,\,\text{see \eqref{af}}},
\end{align*}
and this shows that $D_kH^k_{\alpha\beta}=O(2)$. Now let us consider the equation $D_kH^k_{\alpha\beta}=O(2)$ in more detail, where we again use \eqref{af} to write
\begin{align}
O(2)&=D_kH^k_{\alpha\beta}
=\underbrace{\partial_kH^k_{\alpha\beta}+
u^\gamma_k\partial_\gamma H^k_{\alpha\beta}+
u^\gamma_{ki}\partial^i_\gamma H^k_{\alpha\beta}}_{=O(2)}+
u^\gamma_{kij}\underbrace{\partial^{ij}_\gamma H^k_{\alpha\beta}}_{=O(2)}.\label{even}
\end{align}
From \eqref{even} we get that $\partial^{(ij}_\gamma H^{k)}_{\alpha\beta}=0$, where $(ijk)$ means symmetrization in $ijk$, since third order terms must vanish separately. When we set $i=j=k$ then symmetrization in $\partial^{(ij}_\gamma H^{k)}_{\alpha\beta}=0$ reduces to a single term and we get
\begin{empheq}[box=\fbox]{align}
\partial^{ii}_\gamma H^i_{\alpha\beta}=0.\label{even2}
\end{empheq}
Even if we could use the stronger condition $\partial^{(ij}_\gamma H^{k)}_{\alpha\beta}=0$, instead of \eqref{even2}, it turns out that \eqref{even2} is actually sufficient to complete the proof.\\ 
\quad\\
\underline{Step 6 (second order in I)):} With the help of equation I), \eqref{60}, \eqref{af} and \eqref{even2} we show that
\begin{empheq}[box=\fbox]{align}
H^i_{\alpha\beta}=0.\label{61}
\end{empheq}
This step is quite difficult and we therefore prove it separately in Section 6.\footnote{Notice that this step can maybe also be proven with a different method, when using the so-called $d$-fold operator in \cite[p.379]{bib14} or \cite[p.12]{bib16}.}
\newline\\
\underline{Step 7 (second order in II)):} With the help of equation II), \eqref{af} and \eqref{61} we get that
\begin{empheq}[box=\fbox]{align}
H_{\alpha\beta}=0,\nonumber
\end{empheq}
and therefore all Helmholtz conditions are satisfied. \hfill$\square$
\newline\\
It remains to prove Step 3 and Step 6 in the following sections.

\section{Proof of Step 3}
\label{sec:4}
We briefly explain the notation in this section before we start with the proof. From Step 2 in Section 4 we know that $(f_\alpha)$ must be polynomials of degree $\le n$ in second order coordinates, see \eqref{60}. 
Now recall the Helmholtz expression $H^{ij}_{\alpha\beta}=\partial^{ij}_\beta f_{\alpha}-\partial^{ij}_\alpha f_{\beta}$. When $f_\alpha$ and $f_\beta$ have the form \eqref{60}, then $H^{ij}_{\alpha\beta}$ is of degree $\le n-1$ in the second order coordinates and we get
\begin{empheq}[box=\fbox]{align}
\partial^{i_1j_1}_{\gamma_1}\partial^{i_2j_2}_{\gamma_2}...\partial^{i_nj_n}_{\gamma_n}H^{ij}_{\alpha\beta}=0.\label{644}
\end{empheq}
To determine $H^i_{\alpha\beta}$ in more detail, let us introduce the following short notation. We write
\begin{align}
f_\alpha=A^0_\alpha+A^1_\alpha u_{(2)}+A^2_\alpha u_{(2)}u_{(2)}+...+
A^n_\alpha \underbrace{u_{(2)}u_{(2)}...u_{(2)}}_{\text{$n$-times}},\nonumber
\end{align}
where the coefficients $A^{l}_{\alpha}$ are first order functions, and $u^\gamma_{ij}$ will be identified with $u_{(2)}$, that is, we suppress some of the indices in comparison with \eqref{60}. With this short notation we do a symbolic calculation of the following form (the main focus is on the last expression $2D_k\partial^{ik}_\alpha f_{\beta}$)
\begin{align}
H^i_{\alpha\beta}&=\partial^i_\beta f_{\alpha}+\partial^i_\alpha f_{\beta}-2D_k\partial^{ik}_\alpha f_{\beta}=\nonumber\\
&=O_n(2)-2D_k
\partial^{ik}_\alpha (A^0_\beta+A^1_\beta u_{(2)}+A^2_\beta u_{(2)}u_{(2)}+...+A^n_\beta \underbrace{u_{(2)}u_{(2)}...u_{(2)}}_{\text{$n$-times}})=\nonumber\\
&=O_n(2)+D_k(B^{0,k}_{\alpha\beta}+B^{1,k}_{\alpha\beta} u_{(2)}+B^{2,k}_{\alpha\beta} u_{(2)}u_{(2)}+...+B^{n-1,k}_{\alpha\beta} \underbrace{u_{(2)}u_{(2)}...u_{(2)}}_{\text{$(n-1)$-times}})=\nonumber\\
&=O_n(2)+C^0_{\alpha\beta}u_{(3)}+C^1_{\alpha\beta}u_{(2)}u_{(3)}+...+C^{n-2}_{\alpha\beta}\underbrace{u_{(2)}u_{(2)}...u_{(2)}}_{\text{$(n-2)$-times}}u_{(3)},\label{insb}
\end{align}
where $B^{l,k}_{\alpha\beta}$ and $C^l_{\alpha\beta}$ are again first order functions and $u^\gamma_{ijk}$ is identified with $u_{(3)}$. In exact notation, we can write
\begin{empheq}[box=\fbox]{align}
\partial^{i_1j_1}_{\gamma_1}\partial^{i_2j_2}_{\gamma_2}...\partial^{i_{n-1}j_{n-1}}_{\gamma_{n-1}}H^i_{\alpha\beta}=O_1(2)=O(2).\label{bla2}
\end{empheq}
Next, we want to apply second order partial derivatives $\partial^{i_1j_1}_{\gamma_1}\partial^{i_2j_2}_{\gamma_2}...\partial^{i_lj_l}_{\gamma_l}$, $1\le l\le n$, to equation I) and use the conditions \eqref{644} and \eqref{bla2} to derive further restrictions for $H^i_{\alpha\beta}$ and $H^{ij}_{\alpha\beta}$. But, we do not apply all kinds of second order partial derivatives $\partial^{i_1j_1}_{\gamma_1}\partial^{i_2j_2}_{\gamma_2}...\partial^{i_lj_l}_{\gamma_l}$, we only apply those which can be written as
\begin{align*}
&a)\quad (\partial^{JJ}_\gamma)^r:=\underbrace{\partial^{JJ}_{\gamma_1}\partial^{JJ}_{\gamma_2}...\partial^{JJ}_{\gamma_r}}_{\text{$r$-times}},\quad 1\le r\le n\\
&b)\quad \partial^{jJ}_{\delta}(\partial^{JJ}_{\gamma})^{r-1}:=
\partial^{jJ}_\delta\underbrace{\partial^{JJ}_{\gamma_1}...\partial^{JJ}_{\gamma_{r-1}}}_{\text{$(r-1)$-times}},\quad 1\le r\le n-1\\
&c)\quad \partial^{jj}_{\delta}(\partial^{JJ}_\gamma)^{r-1}:=
\partial^{jj}_{\delta}\underbrace{\partial^{JJ}_{\gamma_1}...\partial^{JJ}_{\gamma_{r-1}}}_{\text{$(r-1)$-times}},\quad 1\le r\le n-1
\end{align*}
where $J$ is a placeholder for (in general different) numbers in $\{1,2,...,n\}$ such that $J\ne j$. Here is the reason why we choose these operators. Recall equation I) which is
\begin{align}
\text{I)}\quad (O^\beta_{jk}(1)-u^\beta_{jk})H^k_{\alpha\beta}+(O^\beta_{jkl}(2)-u^\beta_{jkl})H^{kl}_{\alpha\beta}=0,\label{63}
\end{align}
for $j=1,2,...,n$ and $\alpha=1,2,...,m$. The operators in a), b), c) are defined in such a way that they commute with $u^\beta_{jk}$ in \eqref{63}, except this does not hold for $\partial^{jJ}_\delta$ and $\partial^{jj}_\delta$. Notice that there is summation over $k$, $kl$ and $\beta$ in \eqref{63}, but there is no summation over the indice $j$ and this is crucial, that is, this allows us to define operators which (almost all) commute with $u^\beta_{jk}$. Again, the $J$\rq s can take all possible values $J\in\{1,2,...,n\}$ such that $J\ne j$ and this holds for every single indice $J$ individually, that is, the (in general different) $J$\rq s can have different values. For example, for $a)$ we can equivalently write
\begin{align}
(\partial^{JJ}_\gamma)^r=\partial^{J_1J_2}_{\gamma_1}\partial^{J_3J_4}_{\gamma_2}...\partial^{J_{2r-1}J_{2r}}_{\gamma_r},\quad \text{for $J_s\in \{1,2,...,n\}$ and $J_s\ne j$},\label{75}
\end{align}
and a similar notation holds for $b)$ and $c)$. To indicate that $J$ is in general not a fixed number, we will also write $J\text{\rq s}\in \{1,2,...,n\}$ instead of $J\in\{1,2,...,n\}$. We also use the notation $(\gamma)^r=\gamma_1...\gamma_r$, that is, multi-index notation in $(\gamma)^r$ is assumed. Notice that \eqref{75} actually defines a set of operators for certain $J\text{\rq s}\in\{1,2,...,n\}$, $J\text{\rq s}\ne j$ and there are many combination of such $J\text{\rq s}$ which lead to an operator $(\partial^{JJ}_\gamma)^r$. Therefore, let us now clarify the notation for these combinations and how to handle possible further restrictions for the $J\text{\rq s}$.\\
\hspace*{5mm} Below we consider the set of operators $(\partial^{JJ}_\gamma)^r$, $J\text{\rq s}\in\{1,2,...,n\}$, $j\ne J$, but we will also have further restrictions (similar for $\partial^{jJ}_\delta(\partial^{JJ}_\gamma)^{r-1}$ and $\partial^{JJ}_\delta(\partial^{JJ}_\gamma)^{r-1}$). For example, $J\text{\rq s}\in\{2,3,...,n\}$ and $J\text{\rq s}\ne j$. In general, we denote $S_r\subset\{1,2,...,n\}$ as a subset of $\{1,2,...,n\}$ with $r$ elements (we also define $S_0:=\emptyset$ as the empty set). When we have further restrictions, we simply write $J\text{\rq s}\in S_r$, and we do not always add $J\text{\rq s}\ne j$, since this is always be assumed. Notice that $j$ is allowed to be in the set $S_r$, but it can also be excluded. When we write $J\text{\rq s}\in S_{n}$, then there is only one possibility to construct a set $S_n$ with $n$ elements, namely $S_n=\{1,2,...,n\}$. However, when, for example, we write $J\text{\rq s}\in S_{n-1}$, then there are $n$ possibilities to construct a set $S_{n-1}$ with $n-1$ numbers and we also need a notation for all such possibilities. Therefore, let us define $\mathscr{S}_{n-1}$ as the set which consists of all subsets of $\{1,2,...,n\}$ with $n-1$ elements, that is,
\begin{align*}
\mathscr{S}_{n-1}:=\{&\{2,3,...,n\},\{1,3,4...,n\},\{1,2,4,5,...,n\},...,
\{1,2,...,n-1\}\}
\end{align*}
and a similar definition holds for all the other sets $\mathscr{S}_l$, $0\le l\le n$, which consist of subsets of $\{1,2,...,n\}$ with $l$ elements ($\mathscr{S}_0=\emptyset$). This notation allows us to write a large amount of equations in a structured and compact form, for example, see \eqref{laz}, \eqref{laz2}, \eqref{laz3}.\\
\hspace*{0.5cm} Since we now have clarified our notation, we can start with the main part of the proof of Step 3. The proof is based on a sort of induction. We will write Step $3.k$ for the $k$-th step in the induction.
\newline\\
\underline{Step 3.0 (Start of the proof of Step 3):} The starting point of the induction are equations \eqref{644} and \eqref{bla2}. To get the induction going, we actually only need the weaker conditions written in the following boxes
\begin{empheq}[box=\fbox]{align}
a)&\quad(\partial^{JJ}_\gamma)^{n}H^{kl}_{\alpha\beta}=0,\quad\quad\quad\;\; J\text{\rq s}\in S_n,\quad\quad S_n\in\mathscr{S}_n,\quad\quad\; j=1,2,...,n,\label{laz}\\
b)&\quad(\partial^{JJ}_\gamma)^{n-1}H^{J}_{\alpha\beta}=O(2),\quad J\text{\rq s}\in S_{n-1},\quad S_{n-1}\in\mathscr{S}_{n-1},\; j=1,2,...,n,\label{laz2}\\
c)&\quad(\partial^{JJ}_\gamma)^{n-1}H^{j}_{\alpha\beta}=O(2),\quad J\text{\rq s}\in S_{n-1},\quad S_{n-1}\in\mathscr{S}_{n-1},\; j=1,2,...,n.\label{laz3}
\end{empheq}
By convention, \eqref{laz}, \eqref{laz2} and \eqref{laz3} hold for all $k,l=1,2,...,n$, for all $\alpha,\beta=1,2,...,m$ and for all multi-indices $\gamma$, but we do not explicitly notice it because of lack of space (also in the formulas below). Also notice that the condition $J\text{\rq s}\in S_n$, $S_n\in\mathscr{S}_n$ and the condition $J\text{\rq s}\in S_{n-1}$, $S_{n-1}\in\mathscr{S}_{n-1}$ are equivalent in this case, since $J$ cannot take $n$ different values, because we always assume $J\ne j$ for some $j$. We choose this notation here to keep the same structure in the induction later on.\\
\quad\\
\underline{Step 3.1:} We apply the operators $(\partial^{JJ}_\gamma)^{n-1}$, $J\text{\rq s}\in S_{n-1}$, $S_{n-1}\in\mathscr{S}_{n-1}$ to
equation \eqref{63}, and since we can commute them with $O(1)$, $u^\beta_{jk}$ and $u^\beta_{jkl}$, we get
\begin{align}
a)\quad 0=(O(1)-u^\beta_{jk})\underbrace{(\partial^{JJ}_\gamma)^{n-1}H^k_{\alpha\beta}}_{A}+
\underbrace{(\partial^{JJ}_\gamma)^{n-1}[O(2)H^{kl}_{\alpha\beta}]}_{B}
-u^\beta_{jkl}\underbrace{(\partial^{JJ}_\gamma)^{n-1}H^{kl}_{\alpha\beta}}_{C}.\nonumber
\end{align}
The term $A$ is $O(2)$, see \eqref{laz2} and \eqref{laz3}. The term $B$ is also $O(2)$, since $H^{kl}_{\alpha\beta}$ is $O(2)$. Moreover, the term $C$ is $O(2)$, and since it is the only term which contains a third order coordinate $u^\beta_{jkl}$, which can be varied independently of all other coordinates on $J^kE$, and $H^{kl}_{\alpha\beta}$ is symmetric in $k,l$, we get 
\begin{align*}
\text{Result a):}\quad(\partial^{JJ}_{\gamma})^{n-1}H^{kl}_{\alpha\beta}=0,\quad J\text{\rq s}\in S_{n-1},\quad S_{n-1}\in\mathscr{S}_{n-1},\quad j=1,2,...,n,
\end{align*}
which is the result of Step 3.1 a).\footnote{For example, for $n=2$ we get $\partial^{11}_\gamma H^{kl}_{\alpha\beta}=0$ and $\partial^{22}_\gamma H^{kl}_{\alpha\beta}=0$, but we do not get $\partial^{12}_\gamma H^{kl}_{\alpha\beta}=0$, since we can only choose from a set $S_{n-1}$ of $n-1=1$ numbers.}\\
\hspace*{0.5cm} Then we apply the operators $\partial^{jJ}_\delta(\partial^{JJ}_\gamma)^{n-2}$, $J\text{\rq s}\in S_{n-2},\; S_{n-2}\in\mathscr{S}_{n-2}$ to equation \eqref{63} and we get
\begin{align}
b)\quad 0=&O(1)\underbrace{\partial^{jJ}_\delta(\partial^{JJ}_\gamma)^{n-2}H^k_{\alpha\beta}}_{A}-
\partial^{jJ}_\delta[u^\beta_{jk}(\partial^{JJ}_\gamma)^{n-2}H^k_{\alpha\beta}]+\nonumber\\
&+\underbrace{\partial^{jJ}_\delta(\partial^{JJ}_\gamma)^{n-2}[O(2)H^{kl}_{\alpha\beta}]}_{B}-
u^\beta_{jkl}\underbrace{\partial^{jJ}_\delta(\partial^{JJ}_\gamma)^{n-2}H^{kl}_{\alpha\beta}}_{C}=\nonumber\\
=&O(2)-(\partial^{JJ}_\gamma)^{n-2}H^J_{\alpha\delta}-
u^\beta_{jk}\underbrace{\partial^{jJ}_\delta(\partial^{JJ}_\gamma)^{n-2}H^k_{\alpha\beta}}_{D}=O(2)-(\partial^{JJ}_\gamma)^{n-2}H^J_{\alpha\delta}.\label{676}
\end{align}
The term $A$ and $D$ is $O(2)$, see \eqref{laz2} and \eqref{laz3}. The term $B$ is also $O(2)$ since $H^{kl}_{\alpha\beta}$ is $O(2)$. The term $C$ vanishes, since when we choose the $J\text{\rq s}$ in a set $S_{n-2}$, then the set $j\cup S_{n-2}$ is a set $S_{n-1}$ and we can apply the result of Step 3.1 a). Therefore, we get 
\begin{align*}
\text{Result b):}\quad(\partial^{JJ}_\gamma)^{n-2}H^J_{\alpha\delta}=O(2),\quad J\text{\rq s}\in S_{n-2},\quad S_{n-2}\in\mathscr{S}_{n-2},\quad j=1,2,...,n,
\end{align*}
which is the result of Step 3.1 b).\\
\hspace*{0.5cm} Next, we apply the operators $\partial^{jj}_\delta(\partial^{JJ}_\gamma)^{n-2}$, $J\text{\rq s}\in S_{n-2}$, $S_{n-2}\in\mathscr{S}_{n-2}$ to equation \eqref{63} and we get
\begin{align*}
c)\quad 0=&O(1)\underbrace{\partial^{jj}_\delta(\partial^{JJ}_\gamma)^{n-2}H^k_{\alpha\beta}}_{A}-
\partial^{jj}_\delta[u^\beta_{jk}(\partial^{JJ}_\gamma)^{n-2}H^k_{\alpha\beta}]+\nonumber\\
&+\underbrace{\partial^{jj}_\delta(\partial^{JJ}_\gamma)^{n-2}[O(2)H^{kl}_{\alpha\beta}]}_{B}-
u^\beta_{jkl}\underbrace{\partial^{jj}_\delta(\partial^{JJ}_\gamma)^{n-2}H^{kl}_{\alpha\beta}}_{C}=\nonumber\\
=&O(2)-(\partial^{JJ}_\gamma)^{n-2}H^j_{\alpha\delta}-
u^\beta_{jk}\underbrace{\partial^{jj}_\delta(\partial^{JJ}_\gamma)^{n-2}H^k_{\alpha\beta}}_{D}=O(2)-(\partial^{JJ}_\gamma)^{n-2}H^j_{\alpha\delta},\nonumber
\end{align*}
where the terms $A,B,C,D$ have the same properties as already explained in $b)$, see \eqref{676}. Therefore, we get
\begin{align*}
\text{Result c):}\quad(\partial^{JJ}_\gamma)^{n-2}H^j_{\alpha\delta}=O(2),\quad J\text{\rq s}\in S_{n-2},\quad S_{n-2}\in\mathscr{S}_{n-2},\quad j=1,2,...,n,
\end{align*}
which is the result of Step 3.1 c).\\
\hspace*{0.5cm} Together, the results in $a)$, $b)$ and $c)$ deliver
\begin{empheq}[box=\fbox]{align}
(\partial^{JJ}_\gamma)^{n-1}H^{kl}_{\alpha\beta}&=0,\quad\quad\;\; J\text{\rq s}\in S_{n-1},\quad S_{n-1}\in\mathscr{S}_{n-1},\quad j=1,2,...,n,\label{erstS2}\\
(\partial^{JJ}_\gamma)^{n-2}H^k_{\alpha\delta}&=O(2),\quad J\text{\rq s}\in S_{n-2},\quad S_{n-2}\in\mathscr{S}_{n-2},\quad j=1,2,...,n.\label{erstS}
\end{empheq}
Notice that we always have to prove $a)$ first, since it is needed to prove $b)$ and $c)$.\\
\quad\\
\underline{Step 3.2:} Now we are ready to do this inductively, that is, we can repeat exactly the same argument in every step, and we always use the result in the box from the previous step to derive the result in the actual step. We briefly write down Step 3.2 once again.\\
\hspace*{0.5cm} In Step 3.2 a), we apply the operators
$(\partial^{JJ}_\gamma)^{n-2}$, $J\text{\rq s}\in S_{n-2}$, $S_{n-2}\in\mathscr{S}_{n-2}$ to equation \eqref{63}, and with the result from the previous step (see \eqref{erstS2} and \eqref{erstS}) we derive $(\partial_{u^\gamma_{JJ}})^{n-2}H^{kl}_{\alpha\beta}=0$, $J\text{\rq s}\in S_{n-2}$, $S_{n-2}\in\mathscr{S}_{n-2}$, $j=1,2,...,n$.\\
\hspace*{0.5cm} Then in Step 3.2 $b)$, we apply $\partial^{jJ}_\delta(\partial^{JJ}_\gamma)^{n-3}$, $J\text{\rq s}\in S_{n-3}$, $S_{n-3}\in\mathscr{S}_{n-3}$ to equation \eqref{63} which leads to $(\partial^{JJ}_\gamma)^{n-3}H^J_{\alpha\delta}=O(2)$, $J\text{\rq s}\in S_{n-3}$, $S_{n-3}\in\mathscr{S}_{n-3}$, $j=1,2,...,n$. Again, if we choose the $J\text{\rq s}$ in $S_{n-3}$ then $j\cup S_{n-3}$ can be considered as a set $S_{n-2}$ and we can apply \eqref{erstS} and the result in Step 3.2 a).\\ 
\hspace*{0.5cm} In Step 3.2 $c)$, we apply $\partial^{jj}_\delta(\partial^{JJ}_\gamma)^{n-3}$, $J\text{\rq s}\in S_{n-3}$, $S_{n-3}\in\mathscr{S}_{n-3}$ to equation \eqref{63} and this leads to $(\partial^{JJ}_\gamma)^{n-3}H^j_{\alpha\delta}=O(2)$, $J\text{\rq s}\in S_{n-3}$, $S_{n-3}\in\mathscr{S}_{n-3}$, $j=1,2,...,n$.\\
\hspace*{0.5cm} Together, with the results $a)$, $b)$, $c)$ we get
\begin{empheq}[box=\fbox]{align}
(\partial^{JJ}_\gamma)^{n-2}H^{kl}_{\alpha\beta}&=0,\quad\quad\;\; J\text{\rq s}\in S_{n-2},\quad S_{n-2}\in\mathscr{S}_{n-2},\quad j=1,2,...,n,\nonumber\\
(\partial^{JJ}_\gamma)^{n-3}H^k_{\alpha\delta}&=O(2),\quad J\text{\rq s}\in S_{n-3},\quad S_{n-3}\in\mathscr{S}_{n-3},\quad j=1,2,...,n.\nonumber
\end{empheq}
We repeat these arguments
\begin{align*}
\vdots
\end{align*}
until we get Step $3.(n-2)$.\\
\quad\\
\underline{Step $3.(n-2)$:} In this step we get
\begin{empheq}[box=\fbox]{align}
(\partial^{JJ}_\gamma)^{n-(n-2)}H^{kl}_{\alpha\beta}&=\partial^{JJ}_{\gamma_1}\partial^{JJ}_{\gamma_2} H^{kl}_{\alpha\beta}=0,\quad J\text{\rq s}\in S_{2},\; S_2\in\mathscr{S}_2,\;j=1,...,n,\nonumber\\
(\partial^{JJ}_\gamma)^{n-(n-1)}H^k_{\alpha\delta}&=\partial^{JJ}_\gamma H^k_{\alpha\delta}=O(2),\quad J\text{\rq s} \in S_1,\; S_1\in\mathscr{S}_1,\;j=1,...,n.\label{n-2}
\end{empheq}
For example, for $n=3$ we get $\partial^{11}_\gamma H^k_{\alpha\delta}=0$, $\partial^{22}_\gamma H^k_{\alpha\delta}=0$ and $\partial^{33}_\gamma H^k_{\alpha\delta}=0$.
But we do not get one of the mixed partial derivatives $\partial^{12}_\gamma H^k_{\alpha\delta}=0$, $\partial^{13}_\gamma H^k_{\alpha\delta}=0$
or $\partial^{23}_\gamma H^k_{\alpha\delta}=0$. This will be crucial in the next step, since applying $\partial^{jJ}_\delta$
does not work any longer.\\
\quad\\
\underline{Step 3.$(n-1)$:} The second last step is different in comparison with all the previous steps, since applying the
$\partial^{jJ}_\gamma$-operators does not work any longer, as already mentioned. Formally, in this step, with $a)$, $b)$ and $c)$, we would get
\begin{empheq}[box=\fbox]{align}
&(\partial^{JJ}_\gamma)^{n-(n-1)}H^{kl}_{\alpha\beta}=\partial^{JJ}_\gamma H^{kl}_{\alpha\beta}=0,\quad J\text{\rq s}\in S_1,\;\; S_1\in\mathscr{S}_1,\;\; j=1,2,...,n,\nonumber\\
&(\partial^{JJ}_\gamma)^{n-n}H^k_{\alpha\delta}=H^k_{\alpha\delta}=O(2),\quad\quad\quad\quad J\text{\rq s}\in S_0=\emptyset,\quad j=1,2,...,n,\label{687}
\end{empheq}
and since we cannot choose $J\text{\rq s}\in\emptyset$, there must be something wrong. However, part $a)$ still works, where we get
\begin{align*}
a)\quad 0&=(O(1)-u^\beta_{jk})\underbrace{\partial^{JJ}_\gamma H^k_{\alpha\beta}}_{
\substack{=O(2),\;\text{see\eqref{n-2}}}}+
\underbrace{\partial^{JJ}_\gamma[O(2)H^{kl}_{\alpha\beta}]}_{=O(2)}
-u^\beta_{jkl}\partial^{JJ}_\gamma H^{kl}_{\alpha\beta},
\end{align*}
where $J\text{\rq s}\in S_1$, $S_1\in\mathscr{S}_1$, that is, both of the $JJ$ must be the same now. Therefore, we get 
\begin{align}
\text{Result a):}\quad\partial^{JJ}_\gamma H^{kl}_{\alpha\beta}=0,\quad J\text{\rq s}\in S_1,\quad S_1\in\mathscr{S}_1.\label{diesmalnicht}
\end{align}
\hspace*{0.5cm} Because of \eqref{n-2}, and because of \eqref{diesmalnicht}, Step 3.$(n-1)$ $b)$ does not work any longer, where we would formally get
\begin{align*}
b)\quad 0=&O(1)\underbrace{\partial^{jJ}_\delta H^k_{\alpha\beta}}_{\substack{=O(2),\;\text{only if $j=J$},\\\text{see} \eqref{n-2}}}-
\partial^{jJ}_\delta [u^\beta_{jk}H^k_{\alpha\beta}]+
\underbrace{\partial^{jJ}_\delta [O(2)H^{kl}_{\alpha\beta}]}_{=O(2)}-
u^\beta_{jkl}\underbrace{\partial^{jJ}_\delta H^{kl}_{\alpha\beta}}_{\substack{=0,\;
\text{only if $j=J$},\\\text{see \eqref{diesmalnicht}}}},
\end{align*}
but we always assume $j\ne J$.\\ 
\hspace*{0.5cm} Surprisingly, part $c)$ still works, where we get
\begin{align*}
c)\quad 0=&O(1)\underbrace{\partial^{jj}_\delta H^k_{\alpha\beta}}_{=O(2),\;\text{see} \eqref{n-2}}-
\partial^{jj}_\delta[u^\beta_{jk}H^k_{\alpha\beta}]+
\underbrace{\partial^{jj}_\delta [O(2)H^{kl}_{\alpha\beta}]}_{=O(2)}-
u^\beta_{jkl}\underbrace{\partial^{jj}_\delta H^{kl}_{\alpha\beta}}_{=0,\;\text{see \eqref{diesmalnicht}}}=\\
=&O(2)-H^j_{\alpha\delta}-u^\beta_{jk}\underbrace{\partial^{jj}_\delta H^k_{\alpha\beta}}_{=O(2),\;\text{see \eqref{n-2}}}
\end{align*}
and there does not occur any $J$. Therefore, we get $H^j_{\alpha\delta}=O(2)$.\\
\hspace*{0.5cm} Together with the results in $a)$, $b)$ and $c)$ we get
\begin{empheq}[box=\fbox]{align}
&(\partial^{JJ}_\gamma)^{n-(n-1)}H^{kl}_{\alpha\beta}=\partial^{JJ}_\gamma H^{kl}_{\alpha\beta}=0,\quad J\text{\rq s}\in S_1,\;\; S_1\in\mathscr{S}_1,\;\; j=1,2,...,n,\nonumber\\
&(\partial^{JJ}_\gamma)^{n-n}H^k_{\alpha\delta}=H^k_{\alpha\delta}=O(2),\;\;\quad\quad\text{(formally $J\text{\rq s}\in\emptyset$)},\quad\; j=1,2,...,n.\label{n-1}
\end{empheq}
In the last step, we now also show that $H^{kl}_{\alpha\beta}=0$.
\newline\\
\underline{Step 3.$n$:} Formally, we would apply $(\partial^{JJ}_\gamma)^{n-n}=1$ to equation \eqref{63} and deduce
\begin{align*}
a)\quad 0=&(O(1)-u^\beta_{jk})\underbrace{H^k_{\alpha\beta}}_{=O(2),\,\text{see \eqref{n-1}}}+
\underbrace{O(2)H^{kl}_{\alpha\beta}}_{=O(2)}
-u^\beta_{jkl}H^{kl}_{\alpha\beta},
\end{align*}
and since $H^{kl}_{\alpha\beta}$ is $O(2)$ and symmetric in $k,l$, and third order terms must vanish separately, we get $H^{kl}_{\alpha\beta}=0$, which is the result of Step 3.$n$ a). Equation $b)$ and $c)$ do not provide new information in Step 3.$n$.\\ 
\hspace*{0.5cm} Now \eqref{n-1} and the result of Step 3.$n$ a) provide 
\begin{empheq}[box=\fbox]{align}
H^{kl}_{\alpha\beta}&=0,\nonumber\\
H^k_{\alpha\delta}&=O(2),\nonumber
\end{empheq}
and the proof of Step 3 is complete. \hfill$\square$

\section{Proof of Step 6}
\label{sec:5}
The proof will be again a kind of induction and we will write Step 6.$k$ for the $k$-th step in the induction. Equation I) can be written as
\begin{align}
\text{I)}\quad (O^\beta_{jk}(1)-u^\beta_{jk})H^k_{\alpha\beta}=0,\quad j=1,2,...,n,\quad\alpha=1,2,...,m,\label{66}
\end{align}
when we substitute the expressions \eqref{af} in \eqref{53}. In the following, a partial derivative operator $\partial^{ij}_\alpha$, where $i=j$, is called \textbf{a derivative with same pairs}. Let us consider the differential operator
$\partial^{11}_{\gamma_1}\partial^{22}_{\gamma_2}...\partial^{nn}_{\gamma_n}$ which includes all different kinds of derivatives with same pairs exactly once. Then we define the differential operator
\begin{align*}
\partial^{11}_{\gamma_1}...\wedge^{j_1j_1}...\wedge^{j_2j_2}...\wedge^{j_rj_r}...\partial^{nn}_{\gamma_n},
\end{align*}
where $\wedge^{j_1j_1}$, $\wedge^{j_2j_2}$, ..., $\wedge^{j_rj_r}$ means that the the derivatives $\partial^{j_1j_1}_{\gamma_{j_1}}$, $\partial^{j_2j_2}_{\gamma_{j_2}}$, ..., $\partial^{j_rj_r}_{\gamma_{j_r}}$ are omitted in $\partial^{11}_{\gamma_1}\partial^{22}_{\gamma_2}...\partial^{nn}_{\gamma_n}$. For example,
\begin{align*}
\partial^{11}_{\gamma_1}...\wedge^{jj}...\partial^{nn}_{\gamma_n}=
\partial^{11}_{\gamma_1}\partial^{22}_{\gamma_2}...\partial^{(j-1)(j-1)}_{\gamma_{j-1}}
\partial^{(j+1)(j+1)}_{\gamma_{j+1}}...\partial^{nn}_{\gamma_n}.
\end{align*}
\underline{Step 6.0 (Start of the proof of Step 6):} According to \eqref{even2} we know that
\begin{align}
\boxed{\partial^{jj}_\gamma H^j_{\alpha\beta}=0}\label{allneed}
\end{align}
for all derivatives with same pairs. Notice that there is no summation over $j$ in \eqref{allneed}. Let us consider the expression $\partial^{11}_{\gamma_1}\partial^{22}_{\gamma_2}...\partial^{nn}_{\gamma_n}H^j_{\alpha\beta}$. In the derivatives $\partial^{11}_{\gamma_1}\partial^{22}_{\gamma_2}...\partial^{nn}_{\gamma_n}$, with same pairs, we can always find a $\partial^{jj}_\gamma$-derivative and $\partial^{jj}_\gamma H^j_{\alpha\beta}=0$. Therefore, we get
\begin{empheq}[box=\fbox]{align}
\partial^{11}_{\gamma_1}\partial^{22}_{\gamma_2}...\partial^{nn}_{\gamma_n}H^j_{\alpha\beta}=0.\label{65}
\end{empheq}
Let us now use \eqref{allneed} and \eqref{65} as the starting point in the induction.
\newline\\
\underline{Step 6.1:} We consider equation \eqref{66}, we apply the operator $\partial^{11}_{\gamma_1}\partial^{22}_{\gamma_2}...\partial^{nn}_{\gamma_n}$, and we get (there is no summation over $j$)
\begin{align*}
O^\beta_{jk}(1)\underbrace{\partial^{11}_{\gamma_1}...\partial^{nn}_{\gamma_n}H^k_{\alpha\beta}}_{=0,\;\text{see Step 6.0},\;\eqref{65}}-
\partial^{11}_{\gamma_1}...\wedge^{jj}...\partial^{nn}_{\gamma_n}H^j_{\alpha\gamma_j}-
u^\beta_{jk}\underbrace{\partial^{11}_{\gamma_1}...\partial^{nn}_{\gamma_n}H^k_{\alpha\beta}}_{=0,\;\text{see Step 6.0},\;\eqref{65}}=0.
\end{align*}
which leads to $\partial^{11}_{\gamma_1}...\wedge^{jj}...\partial^{nn}_{\gamma_n}H^j_{\alpha\gamma_j}=0$, where the $\partial^{jj}_{\gamma_j}$-derivative is omitted. For example, for $n=2$, we get $\partial^{22}_\gamma H^1_{\alpha\beta}=0$ and $\partial^{11}_\gamma H^2_{\alpha\beta}=0$, but we do not (yet) get $\partial^{11}_{\gamma}H^1_{\alpha\beta}=0$ and $\partial^{22}_\gamma H^2_{\alpha\beta}=0$. But because of \eqref{allneed}, we also know that $\partial^{11}_{\gamma_1}...\wedge^{kk}...\partial^{nn}_{\gamma_n}H^j_{\alpha\gamma_j}=0$, $k\ne j$, that is, where the $\partial^{jj}_{\gamma_j}$-derivative is included, which together provides
\begin{empheq}[box=\fbox]{align}
\partial^{11}_{\gamma_1}...\wedge^{kk}...\partial^{nn}_{\gamma_n}H^j_{\alpha\gamma_j}=0\quad\text{for all $j,k=1,2,...n$.}\label{einef}
\end{empheq}
Notice that in the result \eqref{einef} of Step 6.1 there is one derivative omitted, in the result of Step 6.2 there will be two derivatives omitted and so on.
\newline\\
\underline{Step 6.2:} Now we can do this inductively, that is, we can repeat exactly the same argument from the previous step. Let us formulate this once again. We consider equation \eqref{66} and we apply the differential operator $\partial^{11}_{\gamma_1}...\wedge^{kk}...\partial^{nn}_{\gamma_n}$, since we want to use the result \eqref{einef} from the previous step to derive further conditions. Thus, we consider the equation
\begin{align}
\text{I)}\quad \partial^{11}_{\gamma_1}...\wedge^{ll}...\partial^{nn}_{\gamma_n}[(O^\beta_{jk}(1)-u^\beta_{jk})H^k_{\alpha\beta}]=0,\label{69}
\end{align}
where $l=1,2,...,n$. If $j=l$ in \eqref{69}, then we can commute all derivatives with $O^\beta_{jk}(1)$ and $u^\beta_{jk}$ and we do not get any new information. Therefore, let $j\ne l$, that is, the $\partial^{jj}_{\gamma_j}$-derivative is included in $\partial^{11}_{\gamma_1}...\wedge^{ll}...\partial^{nn}_{\gamma_n}$. Then we get
\begin{align}
O^\beta_{jk}(1)\underbrace{\partial^{11}_{\gamma_1}...\wedge^{ll}...\partial^{nn}_{\gamma_n}H^k_{\alpha\beta}}
_{=0,\;\text{see Step 6.1,\;\eqref{einef}}}-
\partial^{11}_{\gamma_1}...\wedge^{jj}...&\wedge^{ll}...\partial^{nn}_{\gamma_n}H^j_{\alpha\gamma_j}-\nonumber\\
&-u^\beta_{jk}\underbrace{\partial^{11}_{\gamma_1}...\wedge^{ll}...\partial^{nn}_{\gamma_n}H^k_{\alpha\beta}}
_{\substack{=0,\;\text{see Step 6.1,\;\eqref{einef}}}}=0\label{ueb}
\end{align}
which leads to $\partial^{11}_{\gamma_1}...\wedge^{jj}...\wedge^{ll}...\partial^{nn}_{\gamma_n}H^j_{\alpha\gamma_j}=0$, $j\ne l$, that is, the derivatives $\partial^{jj}_{\gamma_j}$, $\partial^{ll}_{\gamma_l}$ are omitted.\footnote{Notice that in the sum over $k$ in \eqref{ueb} there exists a $k$ such that $k=l$ and therefore it is not possible to apply \eqref{allneed} and we definitely need the result from Step 6.1 in \eqref{einef}.} Again, because of \eqref{allneed}, we can also write $\partial^{11}_{\gamma_1}...\wedge^{rr}...\wedge^{ll}...\partial^{nn}_{\gamma_n}H^j_{\alpha\gamma_j}=0$, $r\ne l$, and $r,l\ne j$, that is, the $\partial^{jj}_{\gamma_j}$-derivative is included in $\partial^{11}_{\gamma_1}...\wedge^{rr}...\wedge^{ll}...\partial^{nn}_{\gamma_n}$. Together, we get
\begin{empheq}[box=\fbox]{align}
\partial^{11}_{\gamma_1}...\wedge^{rr}...\wedge^{ll}...\partial^{nn}_{\gamma_n}H^j_{\alpha\gamma_j}=0,\quad\text{for all $j,r,l=1,2,...,n$,\quad $r\ne l$.}\label{3s}
\end{empheq}
\underline{Step $6.3$:} When we do the procedure we get
\begin{empheq}[box=\fbox]{align}
&\partial^{11}_{\gamma_1}...\wedge^{kk}...\wedge^{rr}...\wedge^{ll}...\partial^{nn}_{\gamma_n}H^j_{\alpha\gamma_j}=0,\quad\text{for all $k,r,l,j=1,2,...,n$,}\nonumber\\
&\quad\quad\quad\quad\quad\quad\quad\quad\quad\quad\quad\quad\quad\quad\quad\quad\quad\;\text{$k,r,l$ are different.}\nonumber
\end{empheq}
We repeat exactly the same argument from the $6.k$-th in the $6.(k+1)$-th Step
\begin{align*}
\vdots
\end{align*}
until we get Step $6.(n-1)$.
\newline\\
\underline{Step $6.(n-1)$:} We do the same calculation as before and we get
\begin{empheq}[box=\fbox]{align}
\partial^{ll}_{\gamma_l}H^j_{\alpha\gamma_j}=0,\quad\text{for all $j,l=1,2,...,n$,}\label{sternstern2}
\end{empheq}
and this is a generalization of \eqref{allneed}, since the equation now also holds for $j\ne l$.
\newline\\
\underline{Step $6.n$:} In the last step we apply the differential operator $\partial^{jj}_{\gamma_j}$ to equation \eqref{66} which leads to
\begin{align*}
O^\beta_{jk}(1)\underbrace{\partial^{jj}_{\gamma_j}H^k_{\alpha\beta}}_{\substack{=0,\,\text{see \eqref{sternstern2}}}}
-H^j_{\alpha\gamma_j}-u^\beta_{jk}
\underbrace{\partial^{jj}_{\gamma_j}H^k_{\alpha\beta}}_{\substack{=0,\,\text{see \eqref{sternstern2}}}}=0,
\end{align*}
and therefore we get
\begin{empheq}[box=\fbox]{align}
H^j_{\alpha\gamma_j}=0,\nonumber
\end{empheq}
which completes the proof of Step 6. \hfill$\square$

\section{Open problems and conclusion}
\label{sec:6}
Concerning Takens\rq{} question, there are several open problems in terms of applications and pure mathematics.\\
\hspace*{0.5cm} When we want to apply Theorem \ref{t47} in the context of physics, as it is motivated in the introduction, we have at least to check two things. First, we have to check whether physically interesting symmetries satisfy the span-condition \eqref{48}. Since the set of symmetries we may assume in physics is usually large, this condition should be satisfied in a lot of interesting cases. Second, we have to check whether the corresponding continuity equations have physical relevance. In both cases, a coordinate invariant definition is in general difficult without further assumptions. For example, energy conservation is usually an equation of the form $u^\alpha_if_\alpha=D_lJ^l_i$, where $(f_\alpha)$ describes a second order differential equation, the characteristic is $(Q^\alpha_i)=(u^\alpha_i)$, and $(J^l_i)$ is the current density, see \eqref{37}. An open question is, what is the coordinate invariant definition of energy conservation, and does an equation of the form $u^\alpha_if_\alpha=D_lJ^l_i$ also make sense for higher order source forms, like fourth order. Also the continuity equations in \cite{bib14,bib16} are only given in local coordinates and the additional assumption $M=\mathbb{R}^n$ is made. In the context of applications, we should also find an explanation why the symmetries $V$ and continuity equations with characteristics $Q$ are connected in the very special form $Q^\alpha=V^\alpha_{ch}$, see Definition \ref{contequ}.\\
\hspace*{0.5cm} Also a very big open question is to find an explanation why a differential equation, given by functions $f_r$, $r=1,2,...,R$, should allow for a source form formulation $\Delta=f_\alpha du^\alpha\wedge dx$, where $R=m$, and where we force a very specific transformation property of $f$, as we explained in Section 3. These are, in a sense, the additional and hidden assumptions in Theorem \ref{t47} (beside the assumptions of symmetries and continuity equations). This also means that we do not consider the symmetries of the function $f$, rather the symmetries of the source form, that is, the weak formulation of the differential equation.\\
\hspace*{0.5cm} Another quite difficult problem, which is also important in the pure mathematical context, is the definition of a variational differential equation in general. When we already assigned a source form to some function $f$ then there is a clear answer whether the source form is variational or not, see Definition \ref{definitionVar}. But when we only consider the differential equation without an assignment to a source form, then there are many equivalent reformulations of the differential equation. Some of them will be variational and some will not. For example, Maxwell\rq s equations are not variational when formulated with the fields $E$ and $B$ as $\nabla E=\rho$, $\nabla B=0$, $\nabla\times E=-\partial_t B$, $\nabla\times B=j+\partial_t E$, but they are variational when formulated with the vector potential $A_\mu$ as $\partial_\mu F^{\mu\nu}=j^\nu$, where $F^{\mu\nu}=\partial^\mu A^\nu-\partial^\nu A^\mu$. Therefore, we would need an explanation why we should consider the one or the other formulation of the differential equation. This problem is also connected to the so-called variational multiplier method which transforms equations equivalently to other differential equations, see \cite{bib28,bib29}. However, the variational multiplier transformation does not cover all equivalent reformulations of differential equations, especially order reduction  or order increase methods are not included, as explained for Maxwell\rq s equations.\\
\hspace*{0.5cm} The question of Takens is very interesting, as motivated in the introduction, and we should probably try to find a reformulation, such that we can avoid the above mentioned problems, or at least some of them. In our opinion, this is the most crucial problem at the time. For example, a reasonable question could be. Is a source form which satisfies certain symmetries and corresponding continuity equations (or some kinds of conservation laws in general) always equivalent to a variational one and can we always assign a source form to any differential expressions $f_r$, $r=1,2,...,R$ in a certain way, such that we get a meaningful weak formulation, for example, by order reduction or order increase methods.

\section{Appendix}
\label{sec:5}
\paragraph{Proof of the conjecture in Section \ref{sec:1}.} We write $u_{(n)}(x)=\frac{d^nu(x)}{dx^n}$. The so-called total derivative operator $\frac{d}{dx}$ is the operator $\frac{d}{dx}:=\frac{\partial}{\partial x}+u_x\frac{\partial}{\partial u}+u_{xx}\frac{\partial}{\partial u_x}+u_{xxx}\frac{\partial}{\partial u_{xx}}+...$, see \eqref{divop} for further details. When we have a conservation law of the form $u_xf=\frac{d}{dx}E$, and when we assume $E=E(x,u,u_x,...,u_{(k)})$ depends on coordinates up to order $k\ge 2$, then the equation
\begin{align}
u_xf=\frac{\partial E}{\partial x}+u_x\frac{\partial E}{\partial u}+u_{xx}\frac{\partial E}{\partial u_{x}}+u_{xxx}\frac{\partial E}{\partial u_{xx}}+...+u_{(k+1)}\frac{\partial E}{\partial u_{(k)}}\label{1}
\end{align}
leads to the following cascade of conditions:
\begin{itemize}
\item The term $u_xf$ on the left hand side in \eqref{1} only depends on the coordinates $(x,u,u_x,u_{xx})$, since we assume that $f$ is of second order. On the right hand side in \eqref{1}, the term $u_{(k+1)}\frac{\partial E}{\partial u_{(k)}}$ must vanish, since otherwise this would be the only term where we have a $u_{(k+1)}$-coordinate, which can be varied independently of the remaining coordinates $(x,u,u_x,...,u_{(k)})$. This forces that $\frac{\partial E}{\partial u_{(k)}}=0$ for all values $(x,u,u_x,...,u_{(k)})$. This means that $E=E(x,u,u_x,...,u_{(k-1)})$.
\item We repeat these arguments until we get $E=E(x,u,u_x)$.
\end{itemize}
Now let us again consider equation \eqref{1}, which now reduces to
\begin{align}
u_xf(x,u,u_x,u_{xx})=\frac{\partial E(x,u,u_x)}{\partial x}+u_x\frac{\partial E(x,u,u_x)}{\partial u}+u_{xx}\frac{\partial E(x,u,u_x)}{\partial u_{x}}.\label{2}
\end{align}
Equation \eqref{2} shows that $f(x,u,u_x,u_{xx})$ must be affine linear in $u_{xx}$, and therefore we can write $f=A+u_{xx}B$ for some functions $A=A(x,u,u_x)$ and $B=B(x,u,u_x)$. More precisely, we get $A=\frac{1}{u_x}(\frac{\partial E}{\partial x}+u_x\frac{\partial E}{\partial u})$ and $B=\frac{1}{u_x}\frac{\partial E}{\partial u_x}$ whenever $u_x\ne 0$.\\
\hspace*{0.5cm} When $f$ is $\frac{\partial}{\partial x}$-invariant, then $f=f(u,u_x,u_{xx})$, and this also means that $A=A(u,u_x)$ and $B=B(u,u_x)$. Therefore, we get $f=A(u,u_x)+u_{xx}B(u,u_x)$. Substituting $f=A+u_{xx}B$ in \eqref{2}, and sorting all terms with respect to the $u_{xx}$-coordinate, and terms which do not involve a $u_{xx}$-coordinate, delivers the two equations
\begin{align}
u_xA(u,u_x)&=\frac{\partial E(x,u,u_x)}{\partial x}+u_x\frac{\partial E(x,u,u_x)}{\partial u},\label{4}\\
u_xB(u,u_x)&=\frac{\partial E(x,u,u_x)}{\partial u_x}.\label{5}
\end{align}
We apply $\frac{\partial}{\partial u_x}$ to equation \eqref{4}, and $\frac{\partial}{\partial x}$, $\frac{\partial}{\partial u}$ to equation \eqref{5}, to deduce
\begin{align}
A+u_x\frac{\partial A}{\partial u_x}&=\frac{\partial^2 E}{\partial x\partial u_x}+\frac{\partial E}{\partial u}+u_x\frac{\partial^2 E}{\partial u\partial u_x},\label{6}\\
0&=\frac{\partial^2 E}{\partial x\partial u_x},\label{7}\\
u_x\frac{\partial B}{\partial u}&=\frac{\partial^2 E}{\partial u\partial u_x}.\label{8}
\end{align}
Substituting \eqref{7} and \eqref{8} in \eqref{6} delivers
\begin{align}
A+u_x\frac{\partial A}{\partial u_x}=\frac{\partial E}{\partial u}+u^2_x\frac{\partial B}{\partial u}.\label{9}
\end{align}
Multiplying equation \eqref{9} by $u_x$, and using \eqref{4}, delivers
\begin{align}
u^2_x\frac{\partial A}{\partial u_x}=-\frac{\partial E}{\partial x}+u^3_x\frac{\partial B}{\partial u}.\label{11}
\end{align}
Applying $\frac{\partial}{\partial x}$ to \eqref{9} delivers
\begin{align}
0=\frac{\partial^2 E}{\partial x\partial u},\label{10}
\end{align}
since $A$ and $B$ do not explicitly depend on $x$. Now equation \eqref{7} and \eqref{10} show that $\frac{\partial E}{\partial x}$ does not depend on $(u,u_x)$, that is, $\frac{\partial E}{\partial x}=\frac{\partial E(x)}{\partial x}$. We use this condition to rewrite equation \eqref{11} as
\begin{align}
u^2_x(\frac{\partial A(u,u_x)}{\partial u_x}-u_x\frac{\partial B(u,u_x)}{\partial u})=-\frac{\partial E(x)}{\partial x}.\label{12}
\end{align}
When we let $u_x\to 0$ in equation \eqref{12}, and when we assume that $\frac{\partial A}{\partial u_x}-u_x\frac{\partial B}{\partial u}$ is a smooth function in $u_x$, then we get $\frac{\partial E(x)}{\partial x}=0$ for all $x$. When $\frac{\partial E}{\partial x}=0$ and when we divide \eqref{2} by $u_x$, then we get
\begin{align}
f(u,u_x,u_{xx})=\frac{\partial E(u,u_x)}{\partial u}+\frac{u_{xx}}{u_x}\frac{\partial E(u,u_x)}{\partial u_x}=\left(\frac{\partial}{\partial u}+\frac{u_{xx}}{u_x}\frac{\partial}{\partial u_x}\right)E,\label{13}
\end{align}
where we know by assumption that $B=\frac{1}{u_x}\frac{\partial E}{\partial u_x}$ is a smooth function (there exists a unique smooth continuation for $u_x\to 0$). We can always write $E$ as $E=L-u_x\frac{\partial L}{\partial u_x}$ for some function $L=L(u,u_x)$, since this differential equation can always be solved for a suitable $L$. Formally, the solution is $L=c(u)u_x-u_x\int\frac{E(u,\tilde u_x)}{\tilde u^2_x}d\tilde u_x$, where $c(u)$ is an arbitrary function depending on $u$ only. That there exists a one time differentiable function $L$ in the whole range of definition is left to the reader. For example, the equation $u_x=L-u_x\frac{\partial L}{\partial u_x}$ has a singular solution $L=-u_x\ln u_x$, where $L$ is not one time differentiable (the resonance case). Notice that $u_xc(u)=\frac{d}{dx}C(u)$ is a trivial Lagrangian, but we can not ignore it, since it is needed to construct a one time differentiable $L$. Now we can write \eqref{13} as 
\begin{align*}
f&=\left(\frac{\partial}{\partial u}+\frac{u_{xx}}{u_x}\frac{\partial}{\partial u_x}\right)\left(L-u_x\frac{\partial L}{\partial u_x}\right)
=\frac{\partial L}{\partial u}-\frac{d}{dx}\frac{\partial L}{\partial u_x},
\end{align*}
which proves that $f$ is variational. When we consider $\frac{dE}{dx}=\frac{dL}{dx}-u_{xx}\frac{\partial L}{\partial u_x}-u_x\frac{d}{dx}\frac{\partial L}{\partial u_x}$ and when we divide through $u_x$, then we also get that $\frac{d}{dx}\frac{\partial L}{\partial u_x}$ must be continues.\hfill$\square$\\
\quad\\
Notice, when we assume that $f$ is $\frac{\partial}{\partial x}$- and $\frac{\partial}{\partial u}$-invariant, then we would not need to make the smoothness assumption of the function $f$ (see Theorem \ref{t47}).

\bibliographystyle{spmpsci}

\bibliography{bibliography/bibSourceForms}

\begin{thebibliography}{10}
\providecommand{\url}[1]{{#1}}
\providecommand{\urlprefix}{URL }
\expandafter\ifx\csname urlstyle\endcsname\relax
  \providecommand{\doi}[1]{DOI~\discretionary{}{}{}#1}\else
  \providecommand{\doi}{DOI~\discretionary{}{}{}\begingroup
  \urlstyle{rm}\Url}\fi

\bibitem{bib7}
Anderson, I.: The Variational Bicomplex.
\newblock Utah State University Technical Report (1989)

\bibitem{bib23}
Anderson, I., Duchamp, T.: {On the Existence of Global Variational Principles}.
\newblock Amer. J. Math. \textbf{102}(5), 781--868 (1980)

\bibitem{bib1}
Anderson, I., Pohjanpelt, J.: {Variational Principles for Natural
  Divergence-free Tensors in Metric Field Theories}.
\newblock J. Geom. Phys. \textbf{62(12)}, 2376--2388 (2012)

\bibitem{bib12}
Anderson, I., Pohjanpelto, J.: {Variational principles for differential
  equations with symmetries and conservation laws I Second order scalar
  equations}.
\newblock Math. Ann. \textbf{299}, 191--222 (1994)

\bibitem{bib13}
Anderson, I., Pohjanpelto, J.: {Variational principles for differential
  equations with symmetries and conservation laws II Polynomial differential
  equations}.
\newblock Math. Ann. \textbf{301}, 627--653 (1995)

\bibitem{bib14}
Anderson, I., Pohjanpelto, J.: {Symmetries, conservation laws and variational
  principles for vector field theories}.
\newblock Math. Proc. Camb. Philos. Soc. \textbf{120}, 369--384 (1996)

\bibitem{bib29}
Anderson, I., Thompson, G.: {The Inverse Problem of the Calculus of Variations
  for Ordinary Differential Equations}.
\newblock Amer. Math. Soc. \textbf{98}(473), 1--110 (1992)

\bibitem{bib52}
Dafinger, M.: {Invariant Source Forms, Conservation Laws, and the Inverse
  Problem of the Calculus of Variations}.
\newblock Dissertation, Carl von Ossietzky University,
  http://oops.uni-oldenburg.de/4067/  (2018)

\bibitem{bib28}
Douglas, J.: {Solution of the Inverse Problem of the Calculus of Variations}.
\newblock Trans. Amer. Math. Soc. \textbf{50}(1), 71--128 (1941)

\bibitem{bib40}
Helmholtz, H.: {Über die physikalische Bedeutung des Prinzips der kleinsten
  Wirkung}.
\newblock Jorn. rein. u. angew. Math. \textbf{100}, 137--166 (1887)

\bibitem{bib31}
Krupka, D.: {Variational sequences in mechanics}.
\newblock Springer-Verlag, Calc. Var. \textbf{5}, 557--583 (1997)

\bibitem{bib36}
Krupka, D.: Intorduction to Global Variational Geomety.
\newblock Atlantis Press (2015)

\bibitem{bib20}
Krupka, D., Saunders, D.: {Jet manifolds and natural bundles}.
\newblock Elsevier B.V., Handbook of Global Analysis pp. 1035--1068 (2008)

\bibitem{bib30}
Krupka, D., \u{S}ed\u{e}nková, J.: {Variational sequences and Lepage forms}.
\newblock Preprint Series in Global Analysis and Applications pp. 1--11 (2004)

\bibitem{bib9}
Krupková, O.: The Geometry of Ordinary Variational Equations.
\newblock Lecture notes in mathematics, 1678, Springer-Verlag Berlin Heidelberg
  New York (1997)

\bibitem{bib18}
Krupková, O., Malíková, R.: {Helmholtz conditions and their
  generalizations}.
\newblock Balkan J. Geom. and Aplications \textbf{15}(1), 80--89 (2010)

\bibitem{bib49}
Lagrange, J.L.: Mécanique Analytique.
\newblock 2nd Edition Paris Courcier 1811-1815 (1788)

\bibitem{bib16}
Manno, G., Pohjanpelto, J., Vitolo, R.: {Gauge invariance, charge conservation,
  and variational principles}.
\newblock J. Geom. Phys. \textbf{58}, 996--1006 (2008)

\bibitem{bib21}
Noether, E.: {Invariante Variationsprobleme}.
\newblock Nachr. Ges. Wiss. Gött. pp. 235--257 (1918)

\bibitem{bib0}
Olver, P.: Applications of Lie Groups to Differential Equations.
\newblock Springer-Verlag New York Berlin Heidelberg, Graduate Texts in
  Mathematics 107 (1986)

\bibitem{bib15}
Pohjanpelto, J.: {Takens\rq{} problem for systems of first order differential
  equations}.
\newblock Ark. Math. \textbf{33}, 343--356 (1995)

\bibitem{bib17}
Takens, F.: {Symmetries, conservation laws and variational principles}.
\newblock Springer-Verlag, New York, Lecture Notes in Mathematics \textbf{597},
  581--603 (1977)

\bibitem{bib45}
Takens, F.: {A Global Version of the Inverse Problem of the Calculus of
  Variations}.
\newblock Jour. Diff. Geom. \textbf{14}, 543--562 (1979)

\bibitem{bib32}
Volná, J., Urban, Z.: {The interior Euler-Lagrange operator in field theory}.
\newblock Preprint Ser. Glob. Var. Geom., Lepage Research Institut pp. 1--14
  (2013)

\end{thebibliography}


\end{document}